\documentclass[a4paper,reqno,11pt]{amsart}
\usepackage[rgb,dvipsnames]{xcolor}
\usepackage[margin=1.15in]{geometry}
\usepackage{amssymb,amsmath,esint,mathrsfs,bm,cite}
\usepackage{graphicx,epstopdf,caption,subcaption,diagbox,multirow,arydshln}
\usepackage{pgfplots}
\pgfplotsset{compat=newest,width=12.0cm}
\usepackage{tikz}
\usetikzlibrary{patterns}
\usepackage{verbatim}
\usepackage[colorlinks=true,linkcolor=MidnightBlue,citecolor=BrickRed,urlcolor=MidnightBlue]{hyperref}
\allowdisplaybreaks[4]
\numberwithin{equation}{section}
\newtheorem{theorem}{Theorem}[section]
\newtheorem{corollary}[theorem]{Corollary}
\newtheorem{lemma}[theorem]{Lemma}
\newtheorem{prop}[theorem]{Proposition}

\theoremstyle{definition}
\newtheorem{remark}[theorem]{Remark}

\newcommand{\RN}[1]{\textup{\uppercase\expandafter{\romannumeral#1}}}

\def\bC{\mathbb{C}} 
\def\bD{\mathbb{D}} 
\def\bE{\mathbb{E}}

\def\bH{\mathbb{H}} \def\cH{\mathcal{H}}
\def\bI{\mathbb{I}} \def\cI{\mathcal{I}}
 
\def\bK{\mathbb{K}} \def\cK{\mathcal{K}}
\def\bL{\mathbb{L}}

\def\bP{\mathbb{P}} 
\def\bQ{\mathbb{Q}} 
\def\bR{\mathbb{R}} 
\def\bS{\mathbb{S}} \def\cS{\mathcal{S}}

\def\bV{\mathbb{V}}

\def\rd{\mathrm{d}}\def\ri{\mathrm{i}}

\newcommand{\tdet}{\operatorname{det}}
\newcommand{\be}{\begin{equation}}
\newcommand{\ee}{\end{equation}}
\newcommand\lb{\label}\newcommand\eq{\eqref}
\newcommand{\pa}{\partial}

\newcommand{\diag}{\operatorname{diag}}
\newcommand{\y}{\right}
\newcommand{\Om}{\Omega}\newcommand{\om}{\omega}

\begin{document}

\title[Subwavelength Resonances with a Size Defect]
{Subwavelength resonances and localization in finite resonator chains with a size defect}

\author[Y.~Gong]{Yihao Gong}
\address[Y.~Gong]{School of Mathematical Sciences, Beijing Normal University, Laboratory of Mathematics and Complex Systems, Ministry of Education, Beijing 100875, China}
\email{yhgong@mail.bnu.edu.cn}

\author[H.~Li]{Haigang Li}
\address[H.~Li]{School of Mathematical Sciences, Beijing Normal University, Laboratory of Mathematics and Complex Systems, Ministry of Education, Beijing 100875, China}
\email{hgli@bnu.edu.cn}

\author[J.~Zhang]{Junhua Zhang}
\address[J.~Zhang]{School of Mathematics and Computer Sciences, Nanchang University, Nanchang 330031, China}
\email{zhangjunhua@ncu.edu.cn}

\begin{abstract}
We study subwavelength resonances of a finite chain of $N\ge2$ close-to-touching spherical resonators in three dimensions with a single radius defect of relative size $\tau\in(0,1]$ and gap width $\varepsilon$. Taking the high-contrast limit before the narrow-gap limit, we prove that the positive generalized eigenvalues of the weighted chain coupling matrix with respect to the volume matrix determine the leading relative resonant frequencies. After eliminating the finite tails, we obtain an exact characteristic determinant and identify a unique critical ratio $\tau_*\in(0,1)$ at which the largest eigenvalue crosses the upper edge of the homogeneous recurrence band $[0,4]$. For every $0<\tau<\tau_*$, exactly one eigenvalue lies above this band, and its eigenvector decays geometrically along the tails. As $\tau\to\tau_*^-$, its recurrence decay length diverges like $(\tau_*-\tau)^{-1/2}$. As $\tau\to0$, exactly one eigenvalue diverges like $\tau^{-2}$ and its eigenvector concentrates at the defect; if the defect is interior, one additional positive eigenvalue has order $\tau$. For suitably normalized relative resonant modes, the sequential limit of the $\varepsilon$-scaled gradient supremum in each gap equals the absolute difference between adjacent reduced amplitudes. For the largest-eigenvalue mode, this yields defect-adjacent gradient concentration as the defect shrinks.
\end{abstract}

\keywords{Subwavelength resonances, close-to-touching resonators, capacitance matrix, size defects, defect-induced localization, gradient concentration}

\subjclass[2020]{35J05, 35B34, 35P25, 35C20}

\maketitle

\section{Introduction}
Controlling waves with structures much smaller than the operating wavelength is a basic capability in subwavelength physics. The classical example is the Minnaert resonance of an air bubble in water~\cite{Min33,AF18}, whose resonant wavelength can far exceed the bubble diameter. Coupling $N$ such resonators hybridizes their individual resonances into $N$ subwavelength modes. At leading order in the contrast, layer-potential methods reduce a finite system to a generalized eigenvalue problem involving the capacitance matrix and the resonator volumes~\cite{FA22,ADH24}. Related reductions are discussed in~\cite{ADH26}. For periodic crystals, Floquet--Bloch theory gives a family of quasiperiodic problems on a unit cell, leading to subwavelength spectral bands and gaps~\cite{AFLLYZ17}. For finite arrays, the analysis must connect two levels of description: the discrete spectral data and the continuous resonant field.

A finite homogeneous chain is periodic in its interior but terminated at two ends, and those ends affect the leading spectrum. A radius defect changes both the volume of the affected resonator and its capacitive couplings to its neighbors, so the response cannot be inferred from either change alone. This distinction matters especially for close-to-touching resonators, whose narrow gaps produce strong coupling and large field gradients. Our aim is to determine the radius ratios that produce an above-band mode, to describe the spectrum as the defect shrinks, and to relate the resulting discrete amplitudes to the gradient in each gap.

Geometric defects have been studied extensively in infinite periodic subwavelength media. Ammari, Fitzpatrick, Hiltunen and Yu~\cite{AFHY18} showed that changing the radius of a single bubble in a two-dimensional periodic bubbly crystal can create a localized mode with frequency inside a subwavelength band gap. Their point-defect analysis covers both dilute and non-dilute regimes. Ammari, Hiltunen and Yu~\cite{AG22} subsequently considered a line of smaller bubbles in a dilute two-dimensional crystal. They proved that the resulting defect band intersects the bulk band gap for arbitrarily small size perturbations and lies entirely within the gap once the perturbation is sufficiently large. The band is not flat, allowing propagation along the defect line while the field remains localized in the transverse direction.

Localization also arises from material perturbations and disorder. Ammari, Davies and Hiltunen~\cite{AA24} related localized modes caused by material-parameter perturbations to the spectrum of a generalized capacitance operator, including a single defect in an infinite periodic structure. For a one-dimensional wave equation with interval resonators, Ammari, Barandun and Uhlmann~\cite{ABU25} studied finite chains with spatial disorder and used a discrete Green function to predict localized modes induced by changes in the wave speeds. Related results include topological edge modes in chains of three-dimensional resonators~\cite{AT20}, interface modes in finite one-dimensional chains~\cite{AL24}, and eigenpair formulas for perturbed block Toeplitz matrices describing non-Hermitian dimer systems~\cite{ABL25}. Here a radius defect changes both the volume weight and the singular narrow-gap couplings in a finite chain.

For close-to-touching spheres in three dimensions, the leading gap contribution to the capacitance matrix grows like $|\log(\varepsilon/R)|$, where $\varepsilon$ is the gap width and $R$ is the resonator radius. Ammari, Davies and Yu~\cite{ADY20} established resonant-frequency separation and gradient blow-up for two close-to-touching spherical resonators. Li and Zhao~\cite{LZ23} extended the analysis to pairs of close-to-touching convex resonators. More recently, Li and Zhang~\cite{LZ25} determined leading resonant-frequency distributions for several configurations of closely spaced spherical resonators, including homogeneous chains, and obtained asymptotics of the resonant fields in narrow gaps. For a homogeneous chain, the leading narrow-gap coupling matrix is tridiagonal with nearest-neighbor couplings. We ask how a radius defect affects the above-band threshold, the shrinking-defect spectrum, and the associated gap fields.

We consider $N\ge2$ spherical resonators in three dimensions. All spheres have radius $R$, except for the $m$th, whose radius is $\tau R$ with $0<\tau\le1$, and every adjacent gap has width $\varepsilon$. Taking first the high-contrast limit $\delta\to0$ at fixed geometry and then the narrow-gap limit $\varepsilon\to0$ with $\tau$ fixed, we obtain a reduced generalized eigenvalue problem for a weighted chain coupling matrix and a diagonal volume matrix. Its positive eigenvalues determine the leading relative resonant frequencies. The leading collective squared frequency is instead proportional to the ratio of the total electrostatic capacity to the total resonator volume. The radius defect enters the reduced problem through the adjacent coupling weight $\zeta(\tau)=\frac{2\tau}{1+\tau}$ and the volume weight $\tau^3$.

As the defect shrinks, both its volume and the couplings to its neighbors vanish, but they do so at different rates. The spectrum therefore contains an unbounded branch and, for an interior defect, an additional positive branch approaching zero. We use complementary Schur reductions to isolate the diverging eigenvalue and control the bounded spectrum, while retaining the response of the finite tails, including endpoint effects, in the threshold equation. For the continuous fields, we prove one-sided mode convergence at fixed geometry and use uniform narrow-gap estimates based on~\cite{LiLY19,Li20}.

We derive an exact characteristic determinant and identify a unique critical ratio $\tau_*\in(0,1)$, determined only by the lengths of the subchains adjacent to the defect. For $0<\tau<\tau_*$, exactly one reduced eigenvalue lies above the homogeneous recurrence band $[0,4]$, and its eigenvector decays geometrically along the finite tails. The crossing through $4$ is transverse, and the recurrence decay length diverges like $(\tau_*-\tau)^{-1/2}$ as $\tau\to\tau_*^-$.

As $\tau\to0$, exactly one eigenvalue diverges like $\tau^{-2}$, with its eigenvector concentrating at the defect. For an interior defect, there is an additional soft branch of order $\tau$; after removal of the diverging branch, the remaining spectrum converges to the combined spectra of the homogeneous subchains. We further identify the sequential scaled gradient limit in each gap with the absolute difference between adjacent reduced amplitudes, yielding defect-adjacent gradient concentration for the largest-eigenvalue mode as $\tau\to0$. At a known defect site, strict monotonicity also implies that the limiting normalized squared frequency of the above-band mode determines $\tau$ uniquely.

\subsection{Mathematical setup}\lb{sec:setup}
We consider a finite chain $D:=\bigcup_{j=1}^ND_j\subset\bR^3$ of $N\ge2$ spherical resonators, as illustrated in Figure~\ref{fig:chain-geometry}. Throughout, we write $x=(x^1,x')\in\bR^3$, where $x'=(x^2,x^3)\in\bR^2$. Fix $R>0$, $0<\tau\le1$ and a defect site $m\in\{1,\ldots,N\}$. Set $R_m:=\tau R$ and $R_j:=R$ for $j\ne m$, and let $D_j:=B_{R_j}(x_j)$, where $x_j=(x_j^1,0,0)$ and $x_1^1<\cdots<x_N^1$. Adjacent resonators are separated by a gap $\varepsilon>0$:
\[
x_{j+1}^1-x_j^1=R_j+R_{j+1}+\varepsilon,\quad j=1,\ldots,N-1.
\]
Throughout, the analysis is at fixed geometry when $N,m,R,\tau$ and all sphere centers $x_j$ are held fixed; then $\varepsilon>0$ is fixed too. When $\tau$ varies at a prescribed gap width, the centers vary so that
all adjacent gaps retain that width.

By reflection symmetry, it suffices to consider sites $1\le m\le(N+1)/2$: $m=1$ is an endpoint defect, and $m\ge2$ is an interior defect.
\begin{figure}[htbp]
\centering
\begin{tikzpicture}[x=1.05cm,y=1.05cm]
\definecolor{chainGray}{RGB}{215,215,215}
\def\rBig{0.44}\def\rSmall{0.22}\def\gap{0.09}
\pgfmathsetmacro{\xNear}{\rBig+\rSmall+\gap}
\pgfmathsetmacro{\xOne}{-5.14}
\pgfmathsetmacro{\xTwo}{-4.17}
\pgfmathsetmacro{\xThree}{-3.20}
\pgfmathsetmacro{\xPrev}{5.62}
\pgfmathsetmacro{\xLast}{6.59}
\foreach \x in {\xOne,\xTwo,\xThree,-\xNear,\xNear,\xPrev,\xLast}{
\shade[ball color=chainGray] (\x,0) circle (\rBig);
\draw[black!20,line width=0.35pt] (\x,0) circle (\rBig);}
\shade[ball color=chainGray] (0,0) circle (\rSmall);
\draw[black!20,line width=0.35pt] (0,0) circle (\rSmall);
\foreach \x in {-2.56,-2.34,-2.12,-1.90,-1.68,-1.46,1.46,1.68,1.90,2.12,2.34,2.56,3.00,3.22,3.44,3.88,4.10,4.32,4.54,4.76,4.98}{
\fill[black!55] (\x,0) circle (1.15pt);}
\pgfmathsetmacro{\xGapA}{\xOne+\rBig}
\pgfmathsetmacro{\xGapB}{\xTwo-\rBig}
\foreach \x in {\xGapA,\xGapB,-\rSmall-\gap,-\rSmall}{
\foreach \y in {-0.24,-0.08,0.08,0.24}{
\draw[black,line width=0.65pt] (\x,{\y-0.055})--(\x,{\y+0.055});}}
\node at ({(\xGapA+\xGapB)/2},-0.51) {$\varepsilon$};
\node at ({-\rSmall-\gap/2},-0.51) {$\varepsilon$};
\node at (\xOne,0.70) {$R$};
\node at (0.1,0.60) {$\tau R$};
\foreach \x/\lab in {\xOne/D_1,\xTwo/D_2,0.1/D_m,\xPrev/D_{N-1},\xLast/D_N}{
\node[font=\small] at (\x,-0.76) {$\lab$};}
\end{tikzpicture}
\caption{\small A chain of close-to-touching resonators with a single radius defect.}
\lb{fig:chain-geometry}
\end{figure}
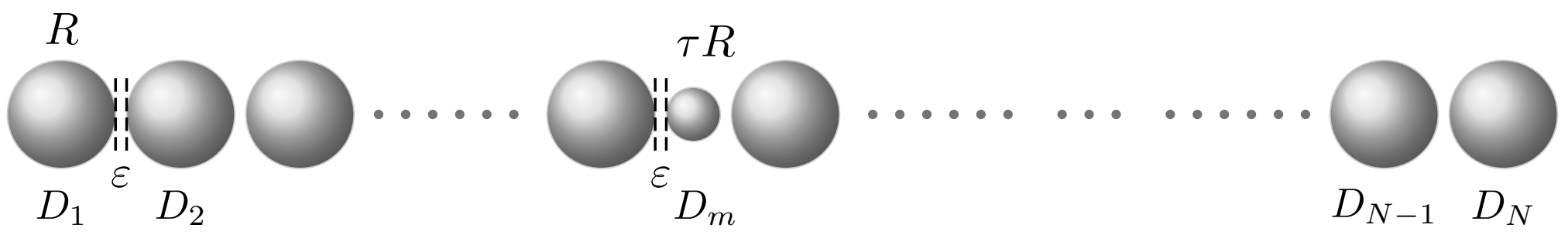
Choose $r_0>0$ and $c>0$ independently of small $\varepsilon$ so that the enlarged contact charts $|x'|<2r_0$ are disjoint and have distance at least $cR$ from every resonator other than the two forming the corresponding gap. For $0<s\le2r_0$, define the narrow region between $D_j$ and $D_{j+1}$ by
\[
\Om_{j,j+1}^{s}:=\{x\in\bR^3\setminus\overline D:x_j^1<x^1<x_{j+1}^1,\ |x'|<s\}.
\]
The constants $r_0/R$ and $c$ may depend on fixed $\tau$. On every compact $\tau$-interval in $(0,1]$, we can choose them uniformly.

Let $\rho,\kappa$ and $\rho_b,\kappa_b$ denote the positive densities and bulk moduli of the background and the resonators, respectively. Define the density contrast and wave speeds by
\[
\delta:=\frac{\rho_b}{\rho},\quad v:=\sqrt{\frac{\kappa}{\rho}},\quad v_b:=\sqrt{\frac{\kappa_b}{\rho_b}}.
\]
The wavenumbers are $k:=\om/v$ and $k_b:=\om/v_b$. Let $\nu$ denote the unit normal on $\pa D$ pointing from $D$ into the exterior domain $\Om:=\bR^3\setminus\overline D$. Subscripts $+$ and $-$ denote the exterior and interior traces, respectively. For an incident field $u^{\mathrm{in}}$ satisfying $(\Delta+k^2)u^{\mathrm{in}}=0$ in $\bR^3$, we write $u=u^{\mathrm{in}}+u^{\mathrm{s}}$ in $\Om$, where $u^{\mathrm{s}}$ is the scattered field. The total field satisfies~\cite{LZ25,LZ23}
\be\lb{trans}
\begin{cases}
\Delta u+k^2u=0&\text{in }\bR^3\setminus\overline D,\\
\Delta u+k_b^2u=0&\text{in }D,\\
u|_+=u|_-&\text{on }\pa D,\\
\delta\pa_\nu u|_+-\pa_\nu u|_-=0&\text{on }\pa D,
\end{cases}
\ee
together with the Sommerfeld radiation condition on the scattered field,
\[
\lim_{r\to\infty}r\big(\pa_r u^{\mathrm{s}}-\ri k u^{\mathrm{s}}\big)=0,\quad r:=|x|,
\]
uniformly in the observation direction for real $k>0$. For complex frequencies, the outgoing condition is understood by analytic continuation of the outgoing fundamental solution.

We study resonances, for which $u^{\mathrm{in}}=0$ and $u=u^{\mathrm{s}}$. A resonant frequency is a nonzero complex $\om$ for which~\eq{trans} admits a nontrivial outgoing solution, called a resonant mode. A branch $\om=\om(\delta)$ is subwavelength if it depends continuously on $\delta$ and $\om(\delta)\to0$ as $\delta\to0$. For each fixed configuration and sufficiently small $\delta>0$, we select the $N$ subwavelength frequencies with nonnegative real parts guaranteed by~\cite[Theorem~2.4]{ADH24}.

We take the high-contrast and close-to-touching limits sequentially,
\be\lb{seq}
\delta\to0\quad\text{for fixed }\varepsilon>0,
\quad\text{followed by}\quad \varepsilon\to0,
\ee
with $N,m,R$, $\tau\in(0,1]$, and the wave speeds $v$ and $v_b$ fixed. We expand the layer potentials as $k\to0$ and use these expansions to follow the resonant branches as $\delta\to0$. Unless stated otherwise, constants and smallness thresholds may depend on the fixed parameters and the Sobolev regularity indices. A subscript $\varepsilon$ indicates that the corresponding constant or convergence rate may depend on $\varepsilon$. After deriving the reduced discrete problem through~\eq{seq} at fixed $\tau>0$, we study $\tau\to0$ in that problem with $N$ and $m$ fixed. For the gap-gradient limits, $\tau$ tends to zero after both limits in~\eq{seq}. The constants in the corresponding shrinking-defect remainder estimates are independent of $\tau$.

For every $n\ge1$, let $\bm1_n,\bm0_n\in\bR^n$ denote the all-ones and zero vectors, respectively. For every $1\le j\le n$, let $\bm e_j\in\bR^n$ denote the $j$th coordinate vector, and let $\bI_n$ denote the identity matrix. We write $\bm1:=\bm1_N$ and drop dimensional subscripts when they are clear from the context. For vectors and matrices, superscripts $T$ and $*$ denote transpose and conjugate transpose, respectively. Unsubscripted norms denote the Euclidean norm for vectors and its induced operator norm for matrices; we also set $\|\bm a\|_\infty:=\max_i|a_i|$.

For every $\bm a=(a_1,\ldots,a_N)^T\in\bR^N$, let $\Phi_{\bm a}$ be the unique exterior harmonic function with trace $a_j$ on $\pa D_j$ and $\Phi_{\bm a}(x)=O(|x|^{-1})$ at infinity. Define $\bC\in\bR^{N\times N}$ by
\[
\bm a^T\bC\bm c:=\int_\Om\nabla\Phi_{\bm a}\cdot\nabla\Phi_{\bm c}\,\rd x,\quad \bm a,\bm c\in\bR^N.
\]
The matrix $\bC$ is symmetric. For complex coefficient vectors, we extend $\Phi_{\bm a}$ by complex linearity; the bilinear identity above remains valid, and
\[
\bm a^*\bC\bm a=\int_\Om|\nabla\Phi_{\bm a}|^2\,\rd x>0
\quad\text{for }\bm a\ne0,
\]
so $\bC$ is positive definite. The total electrostatic capacity is
\[
\operatorname{Cap}(D):=\bm1^T\bC\bm1=\int_\Om|\nabla\Phi_{\bm1}|^2\,\rd x.
\]
The volume matrix $\bV\in\bR^{N\times N}$ is defined by
\be\lb{volmat}
\bV:=\diag\Big(\frac{4\pi}{3}R_1^3,\ldots,\frac{4\pi}{3}R_N^3\Big).
\ee
The capacitance and volume matrices define the generalized eigenvalue problem
\be\lb{capgeig}
\bC\bm a=\lambda\bV\bm a.
\ee
We arrange the eigenvalues of~\eq{capgeig} in nondecreasing order, count them with multiplicity, and use them to label the selected resonant branches. We do not order the complex frequencies themselves.

\subsection{Main results}
Using the coordinate vectors $\bm e_j$ in $\bR^N$, we define the defect-adjacent coupling weight and the dimensionless volume matrix by
\be\lb{dpar}
\zeta(\tau):=\frac{2\tau}{1+\tau},\quad \bD_m^{(\tau)}:=\bI_N+(\tau^3-1)\bm e_m\bm e_m^T.
\ee
The weighted coupling matrix for the chain is
\be\lb{redmat}
\bL_m^{(\tau)}:=\sum_{j=1}^{N-1}s_j^{(\tau)}(\bm e_j-\bm e_{j+1})(\bm e_j-\bm e_{j+1})^T,\quad s_j^{(\tau)}:=\begin{cases}
\zeta(\tau),&m\in\{j,j+1\},\\
1,&\text{otherwise}.
\end{cases}
\ee
Its kernel is $\ker\bL_m^{(\tau)}=\operatorname{span}\{\bm1\}$. The positive eigenvalues of the generalized problem below determine the leading relative resonant frequencies; the collective frequency is governed instead by the bounded part of the capacitance matrix.

\begin{theorem}\lb{thm:defect}
For each fixed $\tau\in(0,1]$, the generalized eigenvalue problem
\be\lb{redgeig}
\bL_m^{(\tau)}\bm a=\alpha\bD_m^{(\tau)}\bm a
\ee
has $N$ simple real eigenvalues:
\be\lb{eigorder}
0=\alpha_1(\tau)<\alpha_2(\tau)<\cdots<\alpha_N(\tau).
\ee
For each $j=2,\ldots,N$, the corresponding resonant frequencies $\om_j$ satisfy
\be\lb{freq}
\lim_{\varepsilon\to0}\lim_{\delta\to0}\frac{4R^2\om_j^2}{3v_b^2\delta|\log(\varepsilon/R)|}=\alpha_j(\tau).
\ee
As $\varepsilon\to0$, the first resonant frequency $\om_1$ satisfies
\be\lb{cmean}
\lim_{\delta\to0}\frac{\om_1^2}{\delta v_b^2}
=\frac{3\operatorname{Cap}(D)}{4\pi R^3(N-1+\tau^3)}
+O\Big(\frac{1}{R^2|\log(\varepsilon/R)|}\Big).
\ee
On every compact subset of $(0,1]$, the implicit error constant is uniform in $\tau$ under uniform bounds on the normalized geometry.
\end{theorem}

We call $\om_1$ the collective branch and $\om_j$, $j=2,\ldots,N$, the relative branches. 

In the homogeneous interior of each tail, the amplitudes satisfy $-a_{\ell-1}+2a_\ell-a_{\ell+1}=\alpha a_\ell$, with characteristic equation $z^2-(2-\alpha)z+1=0$. Its roots are complex conjugates of modulus one for $0<\alpha<4$ and real reciprocal roots for $\alpha>4$. Thus $[0,4]$ is the homogeneous recurrence band; the largest eigenvalue of a finite homogeneous chain is $2(1+\cos(\pi/N))<4$. For $\alpha>4$, define the decaying root and recurrence decay length by
\be\lb{zlocal}
z_-(\alpha):=\frac{2-\alpha+\sqrt{\alpha(\alpha-4)}}{2},\quad \ell_*(\alpha):=\frac{1}{-\log|z_-(\alpha)|},\quad \alpha>4.
\ee
The length $\ell_*$ is measured in chain steps and may exceed a finite tail's length. Let $\vartheta_m$ count the defect neighbors: $\vartheta_m:=1$ at an endpoint and $\vartheta_m:=2$ in the interior.

\begin{theorem}\lb{thm:threshold}
Fix $N\ge2$ and a defect site $m$ as introduced in Subsection~\ref{sec:setup}. Let $L_1,\ldots,L_{\vartheta_m}$ denote the numbers of resonators in the subchains obtained by deleting the defect, with $L_1:=N-1$ at an endpoint and $(L_1,L_2):=(m-1,N-m)$ in the interior. There exists a unique $\tau_*\in(0,1)$ satisfying
\be\lb{crits}
\sum_{i=1}^{\vartheta_m}\frac{2L_i}{2L_i+\tau_*}=2\tau_*^2.
\ee
The largest reduced eigenvalue satisfies
\be\lb{thres}
\alpha_N(\tau)>4\quad(0<\tau<\tau_*),\quad \alpha_N(\tau_*)=4,\quad \alpha_N(\tau)<4\quad(\tau_*<\tau\le1).
\ee
For every $0<\tau<\tau_*$, this eigenvalue is the only one above $4$; the corresponding eigenvector satisfies $a_m\ne0$, and its components decay geometrically along the homogeneous tails. Moreover, $c_*:=-\alpha_N'(\tau_*)>0$, and, as $\tau\to\tau_*^-$,
\[
\ell_*(\alpha_N(\tau))=c_*^{-1/2}(\tau_*-\tau)^{-1/2}\big(1+O(\tau_*-\tau)\big).
\]
\end{theorem}

At an endpoint, the threshold equation reduces to $\tau_*^3+2(N-1)\tau_*^2-(N-1)=0$. For $\tau_*<\tau\le1$, all positive reduced eigenvalues lie in $(0,4)$ and their eigenvector tails are oscillatory; at $\tau=\tau_*$, only the largest eigenvalue reaches $4$. Figure~\ref{fig:defect-spectrum} illustrates the reduced spectra for $N=5$ and $m=1,2,3$, which cover all defect sites up to reflection. Open circles mark the crossings through $4$.

\begin{figure}[t]
\centering
\begingroup
\definecolor{specBlue}{RGB}{31,119,180}
\definecolor{specOrange}{RGB}{255,127,14}
\definecolor{specGreen}{RGB}{44,160,44}
\definecolor{specRed}{RGB}{214,39,40}
\captionsetup[subfigure]{font=footnotesize,skip=2pt}
\pgfplotsset{
specAxis/.style={width=3.75cm,height=2.98cm,scale only axis,
xmin=0,xmax=1,ymin=-.1,ymax=5.25,
xtick={0,.25,.5,.75,1},ytick={0,1,2,3,4,5},
xlabel={$\tau$},ylabel={$\alpha_j(\tau)$},
xlabel style={at={(axis description cs:.5,-.13)},anchor=north},
ylabel style={at={(axis description cs:-.13,.5)},anchor=south},
tick label style={font=\tiny},label style={font=\scriptsize},
ymajorgrids,major grid style={black!9},
axis line style={black!65},tick style={black!55},
clip=true,clip mode=individual,enlargelimits=false,
unbounded coords=jump},
sp2/.style={specBlue,no marks,line width=.95pt},
sp3/.style={specOrange,no marks,line width=.95pt},
sp4/.style={specGreen,no marks,line width=.95pt},
sp5/.style={specRed,no marks,line width=1.3pt}
}
\pgfplotstableread{
t a2 a3 a4 a5 b2 b3 b4 b5 c2 c3 c4 c5
.05 .585755 1.999937 3.414182 762.000125 .061757 1.025258 3.008191 1523.904793 .046486 1.999937 2.048752 1523.904824
.075 .585681 1.999788 3.414107 330.889313 .089245 1.038021 3.012164 661.638349 .067337 1.999788 2.072198 661.638454
.1 .585536 1.999495 3.413959 182.001010 .114708 1.050821 3.016038 363.818433 .086785 1.999497 2.095033 363.818684
.15 .584938 1.998274 3.413335 77.559009 .160182 1.076388 3.023446 154.851095 .121964 1.998293 2.138906 154.851948
.2 .583767 1.995834 3.412062 42.008337 .199290 1.101691 3.030311 83.668708 .152873 1.995919 2.180460 83.670748
.25 .581825 1.991670 3.409821 26.016684 .232949 1.126504 3.036513 51.604034 .180196 1.991937 2.219804 51.608063
.3 .578907 1.985181 3.406175 17.585292 .261916 1.150639 3.041909 34.656647 .204487 1.985870 2.257051 34.663703
.35 .574810 1.975633 3.400486 12.661316 .286804 1.173952 3.046322 24.717411 .226198 1.977187 2.292320 24.728785
.4 .569330 1.962123 3.391767 9.576780 .308113 1.196334 3.049527 18.446026 .245699 1.965291 2.325730 18.463281
.45 .562274 1.943546 3.378413 7.547865 .326244 1.217707 3.051227 14.269020 .263295 1.949513 2.357395 14.293995
.5 .553473 1.918600 3.357677 6.170249 .341523 1.238020 3.051015 11.369441 .279241 1.929110 2.387426 11.404223
.55 .542798 1.885893 3.324747 5.221769 .354214 1.257244 3.048327 9.290628 .293749 1.903280 2.415928 9.337456
.6 .530180 1.844240 3.271657 4.576146 .364537 1.275365 3.042350 7.762193 .307000 1.871206 2.443000 7.823239
.65 .515636 1.793229 3.189038 4.158902 .372678 1.292381 3.031892 6.616658 .319143 1.832137 2.468736 6.693593
.678887732 .506399 1.759824 3.127226 4.000000 .376453 1.301712 3.023145 6.085587 .325705 1.806152 2.483030 6.172010
.7 .499288 1.733897 3.077599 3.913705 .378804 1.308302 3.015179 5.746695 .330307 1.785535 2.493222 5.839916
.75 .481371 1.669029 2.957556 3.780932 .383071 1.323141 2.989576 5.081990 .340604 1.731251 2.516539 5.189384
.8 .462225 1.602576 2.850916 3.709283 .385631 1.336920 2.951357 4.576092 .350127 1.669728 2.538762 4.691383
.85 .442268 1.538443 2.766003 3.668510 .386638 1.349662 2.896011 4.198139 .358958 1.602101 2.559961 4.309430
.884552126 .428241 1.497032 2.719553 3.650274 .386510 1.357876 2.845814 4.000000 .364694 1.552714 2.574046 4.098746
.9 .421954 1.479484 2.701618 3.643858 .386252 1.361395 2.820119 3.926061 .367169 1.530148 2.580200 4.016311
.903211926 .420647 1.475915 2.698090 3.642634 .386184 1.362115 2.814534 3.911740 .367676 1.525428 2.581469 4.000000
.95 .401723 1.427197 2.653649 3.628234 .384638 1.372151 2.724947 3.739870 .374820 1.456034 2.599539 3.791214
1 .381966 1.381966 2.618034 3.618034 .381966 1.381966 2.618034 3.618034 .381966 1.381966 2.618034 3.618034
}\specT
\newcommand{\spRight}{
\node[anchor=west,xshift=2pt,font=\tiny,text=specBlue]at(axis cs:1,.381966){$\frac{3-\sqrt5}{2}$};
\node[anchor=west,xshift=2pt,font=\tiny,text=specOrange]at(axis cs:1,1.381966){$\frac{5-\sqrt5}{2}$};
\node[anchor=west,xshift=2pt,font=\tiny,text=specGreen]at(axis cs:1,2.618034){$\frac{3+\sqrt5}{2}$};
\node[anchor=west,xshift=2pt,font=\tiny,text=specRed]at(axis cs:1,3.618034){$\frac{5+\sqrt5}{2}$};
}
\newcommand{\spGuide}[1]{
\addplot[gray!65,densely dashed,line width=.5pt]
coordinates {(0,4)(1,4)};
\addplot[gray!50,dashdotted,line width=.5pt]
coordinates {(#1,-.1)(#1,5.25)};
\addplot[black!45,line width=.6pt]
coordinates {(0,0)(1,0)};}
\newcommand{\spCurves}[1]{
\addplot[sp2] table[x=t,y=#12]{\specT};
\addplot[sp3] table[x=t,y=#13]{\specT};
\addplot[sp4] table[x=t,y=#14]{\specT};
\addplot[sp5] table[x=t,y=#15]{\specT};
}
\newcommand{\spFinish}[1]{
\addplot[only marks,mark=o,mark size=1.9pt,
mark options={draw=specRed,fill=white}]
coordinates {(#1,4)};
\addplot[only marks,mark=*,mark size=.8pt]
coordinates {(1,.381966)(1,1.381966)(1,2.618034)(1,3.618034)};\spRight}
\begin{subfigure}[t]{.5\textwidth}\centering
\begin{tikzpicture}
\begin{axis}[specAxis,restrict y to domain=-.1:5.25]
\spGuide{.678887732}\spCurves{a}
\addplot[sp2] coordinates {(0,.585786)(.05,.585755)};
\addplot[sp3] coordinates {(0,2)(.05,1.999937)};
\addplot[sp4] coordinates {(0,3.414214)(.05,3.414182)};
\addplot[sp5] coordinates {(.548206042,5.25)(.563206042,5.025511)};
\spFinish{.678887732}
\node[overlay,anchor=east,xshift=-10pt,font=\tiny,text=specBlue]
  at(axis cs:0,.585786){$2-\sqrt2$};
\node[overlay,anchor=east,xshift=-10pt,font=\tiny,text=specGreen]
  at(axis cs:0,3.414214){$2+\sqrt2$};
\end{axis}
\end{tikzpicture}
\caption{Endpoint defect, $m=1$.}
\end{subfigure}\hfill%
\begin{subfigure}[t]{.5\textwidth}\centering
\begin{tikzpicture}
\begin{axis}[specAxis,restrict y to domain=-.1:5.25]
\spGuide{.884552126}\spCurves{b}
\addplot[sp2] coordinates {(0,0)(.05,.061757)};
\addplot[sp3] coordinates {(0,1)(.05,1.025258)};
\addplot[sp4] coordinates {(0,3)(.05,3.008191)};
\addplot[sp5] coordinates {(.736035908,5.25)(.751035908,5.070026)};
\spFinish{.884552126}
\end{axis}
\end{tikzpicture}
\caption{Interior defect, $m=2$.}
\end{subfigure}
\par\vspace{5pt}\noindent
\begin{subfigure}[t]{.5\textwidth}\centering
\begin{tikzpicture}
\begin{axis}[specAxis,restrict y to domain=-.1:5.25]
\spGuide{.903211926}\spCurves{c}
\addplot[sp2] coordinates {(0,0)(.05,.046486)};
\addplot[sp3] coordinates {(0,2)(.05,1.999937)};
\addplot[sp4] coordinates {(0,2)(.05,2.048752)};
\addplot[sp5] coordinates {(.744731458,5.25)(.759731458,5.081804)};
\spFinish{.903211926}
\end{axis}
\end{tikzpicture}
\caption{Centered defect, $m=3$.}
\end{subfigure}
\par\smallskip
{\scriptsize
\tikz\draw[specBlue,line width=.95pt] (0,0)--(.3,0); $\alpha_2$\quad
\tikz\draw[specOrange,line width=.95pt] (0,0)--(.3,0); $\alpha_3$\quad
\tikz\draw[specGreen,line width=.95pt] (0,0)--(.3,0); $\alpha_4$\quad
\tikz\draw[specRed,line width=1.3pt] (0,0)--(.3,0); $\alpha_5$.
}
\endgroup
\caption{Reduced spectra of a chain with different defect sites ($N=5$).}
\lb{fig:defect-spectrum}
\end{figure}
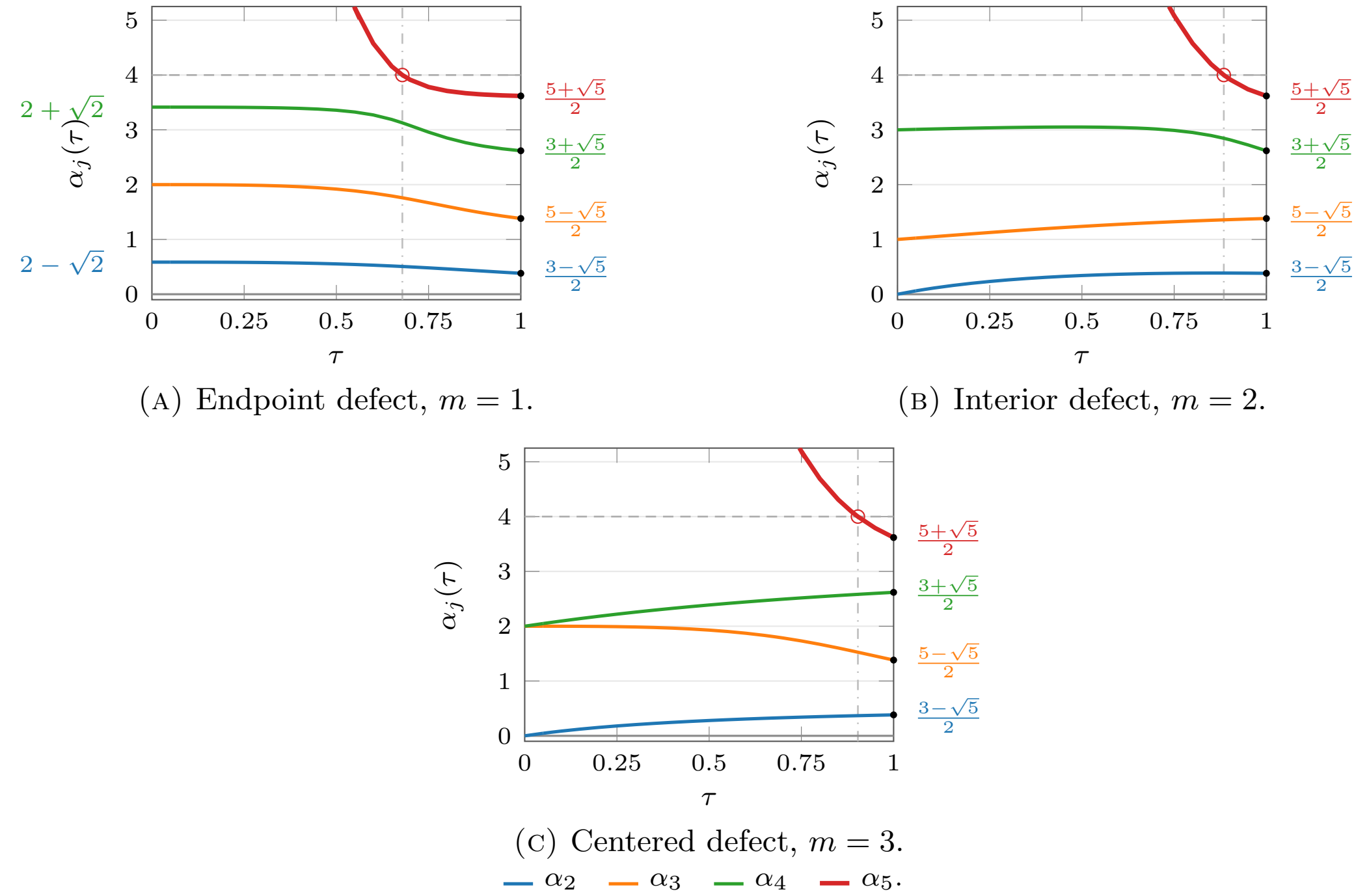

\begin{theorem}\lb{thm:shrink}
Fix $N\ge2$ and a defect site $m$ as in Subsection~\ref{sec:setup}. As $\tau\to0$ in the reduced problem~\eq{redgeig}, exactly one eigenvalue diverges, with
\be\lb{divbr}
\alpha_N(\tau)=\vartheta_m\frac{\zeta(\tau)}{\tau^3}\big(1+O(\tau^3)\big).
\ee
For an interior defect, exactly one additional positive eigenvalue tends to zero, with
\be\lb{softbr}
\alpha_2(\tau)=\frac{N-1}{2(m-1)(N-m)}\zeta(\tau)\big(1+O(\zeta(\tau))\big).
\ee
The remaining $N-1$ eigenvalues, including the soft branch when present, converge with multiplicity to the combined spectrum of the homogeneous subchains left after deleting the defect. Deleting an endpoint leaves one subchain of length $N-1$ and a simple limiting zero eigenvalue. Deleting an interior site leaves subchains of lengths $m-1$ and $N-m$ and a limiting zero eigenvalue of multiplicity two.
\end{theorem}

The differences between the reduced amplitudes on neighboring resonators determine the leading field gradient in each narrow gap.

\begin{theorem}\lb{thm:field}
Fix $\tau\in(0,1]$. For each $j=2,\ldots,N$, choose a real eigenvector $\bm a^{(j)}\in\bR^N$ for $\alpha_j(\tau)$, normalized by $(\bm a^{(j)})^T\bD_m^{(\tau)}\bm a^{(j)}=1$. The corresponding resonant mode $u_j$ can be normalized so that
\be\lb{moderec}
\lim_{\varepsilon\to0}\varepsilon\lim_{\delta\to0}\sup_{\Om_{\ell,\ell+1}^{r_0}}|\nabla u_j|=|a_{\ell+1}^{(j)}-a_\ell^{(j)}|,\quad 1\le\ell<N.
\ee
For $0<\tau<\tau_*$, Theorem~\ref{thm:threshold} gives $a_m^{(N)}\ne0$. We can therefore set $\widetilde u_N:=u_N/a_m^{(N)}$. For the gap $(\ell,\ell+1)$, put $q:=m-\ell$ if $\ell<m$ and $q:=\ell-m+1$ if $\ell\ge m$. Then
\be\lb{phystail}
\lim_{\varepsilon\to0}\varepsilon\lim_{\delta\to0}\sup_{\Om_{\ell,\ell+1}^{r_0}}|\nabla\widetilde u_N|\le2|z_-(\alpha_N(\tau))|^{q-1},\quad 1\le\ell<N.
\ee
\end{theorem}

We normalize the modes in Theorem~\ref{thm:field} by~\eq{dnorm}. For the remaining defect sites, reflection gives the corresponding results with the same normalization and sign or phase conventions.

\begin{remark}\lb{rem:inverse}
We first determine the radius ratio when the defect site is known. For each fixed site, Proposition~\ref{prop:monotone} shows that $\alpha_N(\tau)$ decreases strictly from infinity to $4$ on $(0,\tau_*)$. Hence every limiting normalized squared frequency above $4$ in~\eq{freq} determines a unique $\tau$ in this interval. As $\tau\to0$, we can solve~\eq{divbr} for $\tau$ to obtain
\[
\tau=\sqrt{\frac{2\vartheta_m}{\alpha_N(\tau)}}\big(1+O(\alpha_N(\tau)^{-1/2})\big).
\]
If a homogeneous reference is available, the normalized squared-frequency shift converges to $\alpha_N(\tau)-\alpha_N(1)$ under~\eq{seq}. Since Lemma~\ref{lem:hpath} gives $\alpha_N(1)=2(1+\cos(\pi/N))$, we recover $\alpha_N(\tau)$ by adding this value to the limiting shift. At each prescribed small radius ratio, the leading coefficient in~\eq{divbr} distinguishes an endpoint defect from an interior defect, but it does not distinguish different interior sites.

For a sufficiently small defect, we normalize the largest reduced eigenvector by $a_m=1$. For every $j\ne m$, Lemma~\ref{lem:div} gives $a_j=O(\tau^3)$, so the defect is the unique site where the absolute amplitude is maximal. Because the absolute jumps between adjacent amplitudes tend to one next to the defect and to zero elsewhere, we obtain the same limits for the sequential scaled gap gradients from Corollary~\ref{cor:grad}.

For every candidate site, $\alpha_N$ ranges over $(4,\infty)$ on its localized interval, and reflected sites $m$ and $N+1-m$ have the same spectrum. Hence a single limiting normalized squared frequency does not in general determine both the site and the radius. To recover them from a driven scattering response, we would also need prescribed incident fields and observation data.
\end{remark}

Section~S3.4 of the supplementary material compares the homogeneous chain and the two-resonator case with~\cite{LZ25,ADY20}.

Section~\ref{sec:layer} derives the reduced matrix problem, and Section~\ref{sec:match} establishes its characteristic determinant. Section~\ref{sec:shift} analyzes the spectrum, and Section~\ref{sec:rec} recovers the fields and gap gradients. Detailed auxiliary proofs and additional results are given in the supplementary material.

\section{Capacitance reduction and the reduced matrix problem}\lb{sec:layer}
\subsection{Layer potentials and fixed-geometry reduction}
Fix $\varepsilon>0$ and $\tau\in(0,1]$ as in Subsection~\ref{sec:setup}. Set $X:=H^{-1/2}(\pa D)$ and $Y:=H^{1/2}(\pa D)$ over complex scalars. With the outgoing kernel $G^k(x,y):=-e^{\ri k|x-y|}/(4\pi|x-y|)$, define, initially for smooth densities, the single-layer potential and the adjoint Neumann--Poincar\'e operator by
\[
\begin{aligned}
\cS_D^k[\varphi](x)&:=\int_{\pa D}G^k(x,y)\varphi(y)\,\rd\sigma(y),\quad x\in\bR^3\setminus\pa D,\\
\cK_D^{k,*}[\varphi](x)&:=\int_{\pa D}\pa_{\nu_x}G^k(x,y)\varphi(y)\,\rd\sigma(y),\quad x\in\pa D.
\end{aligned}
\]
The single-layer potential is continuous across $\pa D$, and its trace extends to $\cS_D^k:X\to Y$ with the same notation. With $\cI$ the identity, the normal derivatives satisfy~\cite[Section~2A]{ADH24}
\[
\pa_\nu\cS_D^k[\varphi]|_\pm=\Big(\pm\frac12\cI+\cK_D^{k,*}\Big)[\varphi],\quad \cK_D^{k,*}:X\to X.
\]
Here $+$ and $-$ denote exterior and interior traces, respectively. Write $\cS_D:=\cS_D^0$ and $\cK_D^*:=\cK_D^{0,*}$.

Since $\cS_D:X\to Y$ is invertible, set $\psi_i:=\cS_D^{-1}[\chi_{\pa D_i}]$, where the indicator is one on $\pa D_i$ and zero elsewhere. Component integrals of $X$-densities denote duality pairings with these indicators. The energy identity and uniqueness of the exterior Dirichlet problem give~\cite[Definition~2.5 and Lemma~2.11]{ADH24}
\be\lb{statfield}
[\bC]_{ij}=-\int_{\pa D_i}\psi_j\,\rd\sigma,\quad
\Phi_{\bm a}=\sum_{i=1}^N a_i\cS_D[\psi_i]\quad\text{in }\Om.
\ee
A resonant mode has densities $\phi,\psi\in X$, with $u=\cS_D^{k_b}[\phi]$ in $D$ and $u=\cS_D^k[\psi]$ in $\Om$.

\begin{prop}\lb{prop:red}
For the fixed configuration above, let $\lambda>0$ be a simple generalized eigenvalue of~\eq{capgeig}. There is a unique branch among the $N$ selected subwavelength frequencies corresponding to $\lambda$, and
\[
\om=v_b\sqrt{\delta\lambda}+O_\varepsilon(\delta)\quad(\delta\to0),
\]
where the square root is positive.
\end{prop}

Section~S1 of the supplementary material gives the boundary integral system and small-frequency estimates, and proves Proposition~\ref{prop:red}. These fixed-geometry estimates need not be uniform as $\varepsilon\to0$, so we take the limits in the order~\eq{seq}.

\subsection{Capacitance asymptotics in narrow gaps}\lb{sec:gap}

Let the two adjacent spheres have radii $rR$ and $sR$, where $r,s\in[r_-,r_+]$ and $0<r_-\le r_+<\infty$. Choose $r_0\le c r_-R$ so that the facing spherical graphs are defined over $|x'|<2r_0$. We orient the $x^1$-axis from the left sphere to the right sphere and write their separation as
\[
d_{rs}(x'):=h_+(x')-h_-(x')=\varepsilon+rR-\sqrt{r^2R^2-|x'|^2}+sR-\sqrt{s^2R^2-|x'|^2}.
\]
Set
\[
U_{rs}^{r_0}:=\{(x^1,x'):|x'|<r_0,\ h_-(x')<x^1<h_+(x')\},\quad
J_{rs}(\varepsilon):=\int_{|x'|<r_0}\frac{\rd x'}{d_{rs}(x')}.
\]
For constants $A,B$, we interpolate linearly in the $x^1$-direction and define
\[
\bar v_{A,B}(x):=A+(B-A)\frac{x^1-h_-(x')}{d_{rs}(x')}.
\]

\begin{lemma}\lb{lem:neck}
For every $r,s\in[r_-,r_+]$, as $\varepsilon\to0$ we have, uniformly in $r,s$,
\[
J_{rs}(\varepsilon)=\frac{2\pi Rrs}{r+s}|\log(\varepsilon/R)|+O(R).
\]
If $u\in H^1(U_{rs}^{r_0})$ has constant traces $A,B$ on the facing boundaries, then
\be\lb{fiblo}
\int_{U_{rs}^{r_0}}|\nabla u|^2\,\rd x\ge |A-B|^2J_{rs}(\varepsilon).
\ee
Moreover,
\[
\int_{U_{rs}^{r_0}}|\pa_{x^1}\bar v_{A,B}|^2\,\rd x=|A-B|^2J_{rs}(\varepsilon),\quad
\int_{U_{rs}^{r_0}}|\nabla\bar v_{A,B}|^2\,\rd x\le |A-B|^2(J_{rs}(\varepsilon)+CR).
\]
The constants depend only on $r_-,r_+$ and $r_0/R$, and they do not depend on $\varepsilon,r,s$.
\end{lemma}

\begin{proof}[Proof sketch]
Using the expansion $d_{rs}(x')=\varepsilon+(r+s)|x'|^2/(2Rrs)+O(|x'|^4/R^3)$ to evaluate the gap integral gives the logarithmic coefficient with error $CR$. Applying Cauchy--Schwarz along the $x^1$-direction gives the lower bound, and estimating the transverse derivatives of the linear profile gives the upper bound. The complete proof is given in Section~S2.1 of the supplementary material.
\end{proof}

Write $\mathfrak C(\varepsilon):=\pi R\big|\log(\varepsilon/R)\big|$ and $r_j:=R_j/R$ for $1\le j\le N$, so that $r_m=\tau$ and $r_j=1$ for $j\ne m$. By Lemma~\ref{lem:neck}, we have
\be\lb{gapeng}
J_{r_j,r_{j+1}}(\varepsilon)=\mathfrak C(\varepsilon)s_j^{(\tau)}+O(R).
\ee
By~\eq{redmat}, the sum of the linear-profile gap energies for $\bm a\in\bR^N$ is $\mathfrak C(\varepsilon)\bm a^T\bL_m^{(\tau)}\bm a$ up to $O(R\|\bm a\|^2)$.

For each pair $(D_j,D_{j+1})$, we write the facing boundaries over $|x'|<2r_0$ as $x^1=h_j^-(x')$ and $x^1=h_j^+(x')$, and set
\[
d_{j,j+1}:=h_j^+-h_j^-,\quad \bar v_j(x):=\frac{x^1-h_j^-(x')}{d_{j,j+1}(x')}.
\]
\begin{prop}\lb{prop:cap}
For every chain in Subsection~\ref{sec:setup} whose $\tau$ lies in a fixed compact subset of $(0,1]$, as $\varepsilon\to0$ we have
\be\lb{capexp}
\Big\|\bC-\mathfrak C(\varepsilon)\bL_m^{(\tau)}\Big\|\le CR.
\ee
The estimate is uniform for fixed $N,R$ and uniform positive lower bounds on the normalized chart radii and non-neighboring separations.
\end{prop}

\begin{proof}[Proof sketch]
For the lower bound, we restrict the energy to the disjoint gaps. For the upper bound, we extend the linear profiles through the separated collars. The resulting quadratic-form bounds have error at most $CR\|\bm a\|^2$; symmetry then gives~\eq{capexp}. The complete proof is given in Section~S2.2 of the supplementary material.
\end{proof}

\subsection{Convergence to the reduced spectrum}
The symmetric form of~\eq{redgeig} has matrix
\be\lb{Hddef}
\widehat{\bH}(\tau):=(\bD_m^{(\tau)})^{-1/2}\bL_m^{(\tau)}(\bD_m^{(\tau)})^{-1/2}.
\ee
The matrix $\widehat{\bH}(\tau)$ is real symmetric tridiagonal, with strictly negative subdiagonal entries because the volume and adjacent coupling weights are positive; hence it is irreducible. By~\eq{volmat}, $\bV=(4\pi R^3/3)\bD_m^{(\tau)}$. Set $\bH_\varepsilon:=\bV^{-1/2}\bC\bV^{-1/2}$ and
\[
\widehat{\bH}_{\varepsilon}:=\frac{4R^2}{3|\log(\varepsilon/R)|}\bH_\varepsilon,\quad
\mu_{j,\varepsilon}:=\frac{4R^2\lambda_j(\varepsilon)}{3|\log(\varepsilon/R)|},
\]
where $\lambda_j(\varepsilon)$ are the ordered eigenvalues of~\eq{capgeig}.

\begin{corollary}\lb{cor:srate}
For each fixed $\tau\in(0,1]$, as $\varepsilon\to0$, the symmetric matrices $\widehat{\bH}_\varepsilon$ and $\widehat{\bH}(\tau)$ satisfy
\be\lb{Hrate}
\|\widehat{\bH}_{\varepsilon}-\widehat{\bH}(\tau)\|\le C|\log(\varepsilon/R)|^{-1}.
\ee
Consequently, the eigenvalues satisfy
\be\lb{srate}
|\mu_{j,\varepsilon}-\alpha_j|\le\frac{C}{|\log(\varepsilon/R)|},
\quad j=1,\ldots,N.
\ee
For sufficiently small $\varepsilon$, the eigenvalues of~\eq{capgeig} are simple. Choose corresponding real generalized eigenvectors $\bm a^{(j)}$ and $\bm a_j^\varepsilon$ of~\eq{redgeig} and~\eq{capgeig}, respectively, with unit $\bD_m^{(\tau)}$-norm and $(\bm a_j^\varepsilon)^T\bD_m^{(\tau)}\bm a^{(j)}>0$. Then
\[
\|\bm a_j^\varepsilon-\bm a^{(j)}\|\le C_j|\log(\varepsilon/R)|^{-1},\quad j=1,\ldots,N.
\]
The constants are uniform on each compact $\tau$-interval in $(0,1]$ under uniform geometric bounds. The eigenvector constants also depend on the reduced spectral gaps, whose minimum on such an interval is positive.
\end{corollary}

\begin{proof}[Proof sketch]
For fixed $\tau>0$, multiplying the matrix difference in~\eq{capexp} on both sides by $\bV^{-1/2}$ and applying the scaling in $\widehat{\bH}_\varepsilon$ gives~\eq{Hrate}. Weyl's inequality gives~\eq{srate}, and the reduced spectral gaps give the eigenvector estimate. The complete proof is given in Section~S2.3 of the supplementary material.
\end{proof}

\begin{proof}[Proof of Theorem~\ref{thm:defect}]
The reduced eigenvalue assertions follow from Proposition~\ref{prop:spec}. For sufficiently small $\varepsilon$, Corollary~\ref{cor:srate} allows us to apply Proposition~\ref{prop:red} to each fixed-geometry eigenvalue. Thus $\lim_{\delta\to0}\om_j^2/(\delta v_b^2)=\lambda_j(\varepsilon)$ for $1\le j\le N$. For $j\ge2$, multiplying by $4R^2/(3|\log(\varepsilon/R)|)$ and applying~\eq{srate} proves~\eq{freq}.

For the collective branch, the leading coupling matrix annihilates the constant vector. Set
\[
\bm u:=\frac{\bV^{1/2}\bm1}{(\bm1^T\bV\bm1)^{1/2}}.
\]
Choose an orthonormal basis whose first vector is $\bm u$. In these coordinates, denote the diagonal blocks of $\bH_\varepsilon$ by $a=\bm u^T\bH_\varepsilon\bm u$ and $\bH_+$, and the off-diagonal blocks by $\bm b_\varepsilon^T$ and $\bm b_\varepsilon$. By~\eq{capexp}, $\bH_\varepsilon$ differs by $O(R^{-2})$ from its leading term, which vanishes on $\bm u$ and has a positive gap on $\bm u^\perp$. Thus $\bH_+\ge c|\log(\varepsilon/R)|R^{-2}\bI$ and $\|\bm b_\varepsilon\|\le C/R^2$. Since $\operatorname{Cap}(D)=\bm1^T\bC\bm1$, we have
\[
a=\frac{3\operatorname{Cap}(D)}{4\pi R^3(N-1+\tau^3)},
\quad 0<\lambda_1(\varepsilon)\le a\le C/R^2.
\]
For small $\varepsilon$, $\lambda_1(\varepsilon)<\lambda_{\min}(\bH_+)$, so the first component of its eigenvector is nonzero. The Schur complement gives
\[
0\le a-\lambda_1(\varepsilon)
=\bm b_\varepsilon^T(\bH_+-\lambda_1(\varepsilon)\bI)^{-1}\bm b_\varepsilon
\le\frac{C}{R^2|\log(\varepsilon/R)|}.
\]
This proves~\eq{cmean}, uniformly on compact $\tau$-intervals by the reduced spectral gap and the stated geometric bounds. Since $\mu_{j,\varepsilon}\to\alpha_j$ and the reduced eigenvalues are simple, exactly one selected branch is collective and the remaining $N-1$ branches are relative.
\end{proof}

\section{Characteristic determinant of the finite defect chain}\lb{sec:match}

We derive the characteristic determinant by eliminating the subchains with their free endpoints. The reduced problem defined by~\eq{dpar} and~\eq{redmat} extends algebraically to every $\tau>0$. Set $\bL:=\bL_m^{(\tau)}$ and $\bD:=\bD_m^{(\tau)}$. Because $\bD$ is positive definite, we can apply the min--max principle over real linear subspaces $E\subset\bR^N$:
\be\lb{minmax}
\alpha_j=\min_{\dim E=j}\ \max_{\bm a\in E\setminus\{0\}}\frac{\bm a^T\bL\bm a}{\bm a^T\bD\bm a},\quad\mbox{where}~ \bm a^T\bL\bm a=\sum_{\ell=1}^{N-1}s_\ell^{(\tau)}(a_{\ell+1}-a_\ell)^2.
\ee
We first show that every eigenvalue is simple, so that we can follow each branch as $\tau$ varies.

\begin{prop}\lb{prop:spec}
For every $\tau>0$, the generalized eigenvalue problem~\eq{redgeig} has $N$ simple real eigenvalues, ordered as in~\eq{eigorder}. Its zero eigenspace is $\operatorname{span}\{\bm1\}$, and every eigenvector associated with a positive eigenvalue satisfies $\bm1^T\bD_m^{(\tau)}\bm a=0$.
\end{prop}

\begin{proof}
By~\eq{minmax}, $\bL$ is positive semidefinite with kernel $\operatorname{span}\{\bm1\}$. The nonzero subdiagonal entries of the symmetric matrix~\eq{Hddef} determine an eigenvector recursively from its first component. Thus every eigenvalue is simple. Multiplying~\eq{redgeig} by $\bm1^T$ gives $\bm1^T\bD\bm a=0$ for $\alpha>0$.
\end{proof}

For a homogeneous subchain, set $\bL_{1,0}:=(0)$. For $n\ge2$, let $\bm e_j$ be the coordinate vectors in $\bR^n$ and define
\be\lb{hommat}
\bL_{n,0}:=\sum_{j=1}^{n-1}(\bm e_j-\bm e_{j+1})(\bm e_j-\bm e_{j+1})^T.
\ee
The index $0$ denotes the homogeneous chain.

To eliminate a homogeneous subchain, we need its determinant and the inverse entry at the endpoint adjacent to the defect. We compute both from the polynomials $h_s(\alpha)$ and $q_s(\alpha;\tau)$, defined by
\be\lb{cont}
h_0:=1,~~ h_1:=1-\alpha,~~ h_s:=(2-\alpha)h_{s-1}-h_{s-2}~(s\ge2),~~
q_s:=h_s+(\zeta(\tau)-1)h_{s-1}~(s\ge1).
\ee

For $s\ge1$, order a tail from its free endpoint toward the defect. In $\bL-\alpha\bD$, this tail has the block
\[
\bQ_s(\alpha;\tau):=\bL_{s,0}+\zeta(\tau)\bm e_s\bm e_s^T-\alpha\bI_s,\quad \bm e_s\in\bR^s.
\]
We retain the added diagonal term because the tail remains coupled to the defect.

\begin{lemma}\lb{lem:tdet}
For $s\ge1$ and every $\alpha$, $\tdet\bQ_s(\alpha;\tau)=q_s(\alpha;\tau)$. Whenever $q_s(\alpha;\tau)\ne0$,
\[
\bm e_s^T\bQ_s(\alpha;\tau)^{-1}\bm e_s=\frac{h_{s-1}(\alpha)}{q_s(\alpha;\tau)}.
\]
\end{lemma}

\begin{proof}
For $s=1$, the determinant is $\zeta(\tau)-\alpha=q_1$ and the cofactor is $h_0=1$. For $s\ge2$, the leading principal minors through order $s-1$ satisfy~\eq{cont}. Expanding along the last row gives $\tdet\bQ_s=(1+\zeta(\tau)-\alpha)h_{s-1}-h_{s-2}=q_s$. The last diagonal cofactor is $h_{s-1}$, so Cramer's rule gives the inverse entry.
\end{proof}

Let $\widetilde{\bQ}_s$ denote the block with the subchain order reversed. Reversing both rows and columns gives $\tdet\widetilde{\bQ}_s=q_s$ and $\bm e_1^T\widetilde{\bQ}_s^{-1}\bm e_1=h_{s-1}/q_s$ whenever $q_s\ne0$.

Set $x:=1-\alpha/2$ and define the Chebyshev polynomials by $U_{-1}:=0$, $U_0:=1$ and $U_{s+1}(x):=2xU_s(x)-U_{s-1}(x)$. The sequence $U_s(x)-U_{s-1}(x)$ satisfies the initial conditions and recurrence in~\eq{cont}; hence
\be\lb{cheb}
h_s(\alpha)=U_s(x)-U_{s-1}(x).
\ee
We next compute the homogeneous spectrum from~\eq{cheb}; its largest eigenvalue will bound the retained subchain spectra.

\begin{lemma}\lb{lem:hpath}
For $N\ge2$, the homogeneous chain coupling matrix $\bL_{N,0}$ defined in~\eq{hommat} has eigenvalues
\[
2\Big(1-\cos\frac{(j-1)\pi}{N}\Big),\quad j=1,\ldots,N.
\]
In particular, its largest eigenvalue is $2(1+\cos(\pi/N))<4$ and its kernel is spanned by the constant vector.
\end{lemma}

\begin{proof}
The last-row expansion and~\eq{cheb} give
\be\lb{pathdet}
\tdet(\bL_{N,0}-\alpha\bI_N)=(1-\alpha)h_{N-1}-h_{N-2}=-\alpha U_{N-1}(1-\alpha/2).
\ee
This also holds for $N=2$, when the determinant is $-\alpha(2-\alpha)$ (see~\cite[Lemma~3.2 and (3.8)]{LZ25} for $N>2$). The zeros of $U_{N-1}$ give the positive eigenvalues, and the homogeneous quadratic form vanishes only on constants.
\end{proof}

We now apply Lemma~\ref{lem:tdet} to eliminate each homogeneous subchain and compute the full characteristic determinant.

\begin{prop}\lb{prop:chain}
Fix $N\ge2$, a defect site $m$ in the representative range of Subsection~\ref{sec:setup}, and $\tau>0$. For an interior defect, put $L:=m-1$ and $K:=N-m$. The characteristic determinant of~\eq{redgeig} is
\be\lb{matchi}
F_{N,m}(\alpha;\tau):=\tdet\big(\bL_m^{(\tau)}-\alpha\bD_m^{(\tau)}\big) =(2\zeta(\tau)-\alpha\tau^3)q_Lq_K -\zeta(\tau)^2\big(h_{L-1}q_K+h_{K-1}q_L\big).
\ee
For an endpoint defect at $m=1$,
\be\lb{matche}
F_{N,1}(\alpha;\tau):=\tdet\big(\bL_1^{(\tau)}-\alpha\bD_1^{(\tau)}\big) =(\zeta(\tau)-\alpha\tau^3)q_{N-1}-\zeta(\tau)^2h_{N-2}.
\ee
For the other endpoint, we reflect the chain. For every $\alpha$, including each zero of a subchain determinant, these identities remain valid and give the complete finite-chain spectrum.
\end{prop}

\begin{proof}
Let $\bm e_L\in\bR^L$ select the defect end of the left subchain, and let $\bm e_1\in\bR^K$ select that of the right subchain. In the original resonator order,
\[
\bL-\alpha\bD=\begin{pmatrix}
\bQ_L&-\zeta(\tau)\bm e_L&0\\
-\zeta(\tau)\bm e_L^T&2\zeta(\tau)-\alpha\tau^3&-\zeta(\tau)\bm e_1^T\\
0&-\zeta(\tau)\bm e_1&\widetilde{\bQ}_K
\end{pmatrix}.
\]
By Lemma~\ref{lem:tdet}, the tail determinants are $q_L,q_K$, and their defect-adjacent cofactors are $h_{L-1},h_{K-1}$. When $q_Lq_K\ne0$, we can eliminate both tails and obtain~\eq{matchi}.

When $m=1$, the tail block is $\widetilde{\bQ}_{N-1}$, the defect diagonal entry is $\zeta(\tau)-\alpha\tau^3$, and the coupling is $-\zeta(\tau)\bm e_1$, where $\bm e_1\in\bR^{N-1}$. For $q_{N-1}\ne0$, we eliminate this block to obtain~\eq{matche}; $h_0=1$ also covers $N=2$. Since both sides of each identity are polynomials in $\alpha$, the identities extend to every $\alpha$.
\end{proof}

At $\tau=1$, $q_s=h_s$ and the two formulas reduce to~\eq{pathdet}.

\section{Spectral consequences of a radius defect}\lb{sec:shift}

We use the characteristic determinant to locate the localization threshold, describe the finite tails, and analyze the spectrum as the defect shrinks.

Fix $N$ and $m$ throughout this section. For $\tau>0$, $F_{N,m}(\alpha;\tau)$ is real analytic in $(\alpha,\tau)$ by Proposition~\ref{prop:chain}, and every eigenvalue is a simple root by Proposition~\ref{prop:spec}. We can therefore apply the implicit function theorem: near every $\tau>0$, each ordered branch $\alpha_j(\tau)$ is real analytic.

\subsection{The localization threshold}

To decide when $\alpha_N(\tau)$ crosses $4$, we eliminate the coordinates away from the defect. We shall use $h_s(4)=(-1)^s(2s+1)$, which follows from~\eq{cheb}.

Delete the defect coordinate while keeping the remaining resonators in their original order, and write
\[
\bm a_m^-:=(a_1,\ldots,a_{m-1},a_{m+1},\ldots,a_N)^T\in\bR^{N-1}.
\]
Define $\bL_{N,m}^{\mathrm{ho}}\in\bR^{(N-1)\times(N-1)}$ for the remaining homogeneous subchains, and let $\bm e_m^{\mathrm{ad}}\in\bR^{N-1}$ select the neighbors of the defect:
\be\lb{taildata}
\begin{aligned}
m=1:~&\bL_{N,1}^{\mathrm{ho}}:=\bL_{N-1,0},~\bm e_1^{\mathrm{ad}}:=\bm e_1,\\
m\ge2:~&\bL_{N,m}^{\mathrm{ho}}:=\bL_{m-1,0}\oplus\bL_{N-m,0},~\bm e_m^{\mathrm{ad}}:=\bm e_{m-1}+\bm e_m.
\end{aligned}
\ee
All coordinate vectors in~\eq{taildata} belong to $\bR^{N-1}$. For an interior defect, $\bm e_{m-1}$ selects $a_{m-1}$, and $\bm e_m$ selects $a_{m+1}$. Set $\bE_m^{\mathrm{ad}}:=\diag(\bm e_m^{\mathrm{ad}})\in\bR^{(N-1)\times(N-1)}$.

To retain the couplings at these neighbors, set $\bL_{N,m}^{\mathrm{re}}(\tau):=\bL_{N,m}^{\mathrm{ho}}+\zeta(\tau)\bE_m^{\mathrm{ad}}$. In the order $(a_m,\bm a_m^-)\in\bR\oplus\bR^{N-1}$, we can write $\bL_m^{(\tau)}-\alpha\bD_m^{(\tau)}$ as
\be\lb{perm}
\begin{pmatrix}
\vartheta_m\zeta(\tau)-\alpha\tau^3&-\zeta(\tau)(\bm e_m^{\mathrm{ad}})^T\\
-\zeta(\tau)\bm e_m^{\mathrm{ad}}&\bL_{N,m}^{\mathrm{re}}(\tau)-\alpha\bI_{N-1}
\end{pmatrix}.
\ee
Because we reorder both rows and columns, the determinant and generalized spectrum are unchanged. We continue to index components by the original resonators. For $0<\tau\le1$ and $\alpha\ge4$, we shall show that the lower-right block is invertible and use the scalar Schur complement
\[
g(\alpha;\tau):=\vartheta_m\zeta(\tau)-\alpha\tau^3-\zeta(\tau)^2(\bm e_m^{\mathrm{ad}})^T\big(\bL_{N,m}^{\mathrm{re}}(\tau)-\alpha\bI_{N-1}\big)^{-1}\bm e_m^{\mathrm{ad}}.
\]
Once we bound the retained spectrum below $4$, we can determine the threshold from the sign of $g(4;\tau)$.

\begin{prop}\lb{prop:band}
For $0<\tau\le1$, the spectrum of $\bL_{N,m}^{\mathrm{re}}(\tau)$ lies below $4$. The scalar Schur complement $g(\cdot;\tau)$ is strictly decreasing on $[4,\infty)$, and
\[
\alpha_N(\tau)>4\iff g(4;\tau)>0.
\]
If $g(4;\tau)>0$, exactly one eigenvalue lies above $4$. Otherwise, $\alpha_N(\tau)\le4$, with equality precisely when $g(4;\tau)=0$.
\end{prop}

\begin{proof}
Since $0<\zeta(\tau)\le1$, we have $\bL_{N,m}^{\mathrm{re}}(\tau)\le\bL_{N,m}^{\mathrm{re}}(1)$ in the quadratic-form sense. We obtain $\bL_{N,m}^{\mathrm{re}}(1)$ by deleting the defect coordinate from $\bL_{N,0}$, so we can apply Cauchy's interlacing theorem:
\[
\lambda_{\max}\big(\bL_{N,m}^{\mathrm{re}}(\tau)\big)\le\lambda_{\max}\big(\bL_{N,m}^{\mathrm{re}}(1)\big)\le\lambda_{\max}\big(\bL_{N,0}\big)=2\Big(1+\cos\frac{\pi}{N}\Big)<4.
\]
Hence $\bL_{N,m}^{\mathrm{re}}(\tau)-\alpha\bI_{N-1}$ is negative definite for every $\alpha\ge4$, uniformly for $0<\tau\le1$.

We can therefore eliminate the lower-right block in~\eq{perm} and obtain
\be\lb{blockdet}
\tdet\big(\bL_m^{(\tau)}-\alpha\bD_m^{(\tau)}\big)=\tdet\big(\bL_{N,m}^{\mathrm{re}}(\tau)-\alpha\bI_{N-1}\big)\,g(\alpha;\tau),\quad\alpha\ge4.
\ee
The first factor is nonzero on $[4,\infty)$, so the eigenvalues there are precisely the zeros of $g(\cdot;\tau)$. Differentiating gives
\[
\pa_\alpha g(\alpha;\tau)=-\tau^3-\zeta(\tau)^2(\bm e_m^{\mathrm{ad}})^T\big(\bL_{N,m}^{\mathrm{re}}(\tau)-\alpha\bI_{N-1}\big)^{-2}\bm e_m^{\mathrm{ad}}<0.
\]
As $\alpha\to\infty$, $g(\alpha;\tau)=-\alpha\tau^3+\vartheta_m\zeta(\tau)+O(\alpha^{-1})\to-\infty$. Strict monotonicity gives a unique root above $4$ exactly when $g(4;\tau)>0$, and a root at $4$ exactly when $g(4;\tau)=0$.
\end{proof}

\begin{proof}[Proof of the threshold assertions in Theorem~\ref{thm:threshold}]
We must show that $g(4;\tau)$ changes sign exactly once and that $\alpha_N'(\tau_*)<0$. Proposition~\ref{prop:band} then gives the threshold assertions.

\textbf{Step 1: Evaluate $g$ at $\alpha=4$.} After deleting the defect, we have $\vartheta_m$ subchains of lengths $L_i$, as in Theorem~\ref{thm:threshold}. For either orientation, the inverse entry at the defect end equals $h_{L_i-1}(4)/q_{L_i}(4)$ by Lemma~\ref{lem:tdet}. We substitute $h_s(4)=(-1)^s(2s+1)$, with $h_0=1$ when $L_i=1$, to obtain
\[
\frac{h_{L_i-1}(4)}{q_{L_i}(4)}
=-\frac{2L_i-1}{4L_i-(2L_i-1)\zeta(\tau)}.
\]
For $0<\tau\le1$, this denominator is positive because $\zeta(\tau)\le1$. We can therefore substitute the inverse entries into $g$ and use~\eq{dpar} to obtain
\be\lb{phisign}
g(4;\tau)=2\tau H(\tau),\quad H(\tau):=\sum_{i=1}^{\vartheta_m}\frac{2L_i}{2L_i+\tau}-2\tau^2.
\ee
Thus $g(4;\tau_*)=0$ is equivalent to~\eq{crits}. At an endpoint, we set $L_1=N-1$ and multiply by $2L_1+\tau_*$ to obtain the stated cubic.

\textbf{Step 2: Prove uniqueness and compute the derivative.} For $0<\tau\le1$, we have
\[
H'(\tau)=-\sum_{i=1}^{\vartheta_m}\frac{2L_i}{(2L_i+\tau)^2}-4\tau<0.
\]
Since $H(0^+)=\vartheta_m>0$ and $H(1)<\vartheta_m-2\le0$, $H$ has exactly one zero $\tau_*\in(0,1)$. By~\eq{phisign}, $g(4;\tau)>0$ precisely for $0<\tau<\tau_*$. We then apply Proposition~\ref{prop:band} to obtain~\eq{thres} and $\alpha_N(\tau_*)=4$.

When $\tau=\tau_*$, we have $\pa_\tau g(4;\tau_*)=2\tau_*H'(\tau_*)<0$ and $\pa_\alpha g(4;\tau_*)<0$, the second inequality from Proposition~\ref{prop:band}. The retained block remains invertible near $(4,\tau_*)$, where~\eq{blockdet} makes the eigenvalue equation equivalent to $g(\alpha;\tau)=0$. The implicit function theorem identifies the local solution with $\alpha_N(\tau)$ and gives
\[
\alpha_N'(\tau_*)
=-\frac{\pa_\tau g(4;\tau_*)}{\pa_\alpha g(4;\tau_*)}<0.
\]
Since $\alpha_N$ is analytic and $\alpha_N'(\tau_*)<0$, we can expand it as
\be\lb{crossing}
\alpha_N(\tau)-4=c_*(\tau_*-\tau)+O((\tau-\tau_*)^2).
\ee
\end{proof}

For small $\varepsilon$ and $\tau$ near $\tau_*$, the largest normalized capacitance eigenvalue lies on the same side of $4$ as $\alpha_N(\tau)$ whenever $|\tau-\tau_*|\ge K_0|\log(\varepsilon/R)|^{-1}$ for sufficiently large $K_0$; see Section~S3.1 of the supplementary material. Section~S3.2 treats the algebraic extension to $\tau\ge1$.

\subsection{Localization on finite tails}
For $\alpha>4$, the homogeneous recurrence has two real reciprocal roots. We impose the free-endpoint condition on their linear combination and obtain a reflected term on each finite tail.

\begin{prop}\lb{prop:tail}
Let $(\alpha,\bm a)$ solve~\eq{redgeig} with $\alpha>4$. Then the root defined in~\eq{zlocal} satisfies $z_-(\alpha)\in(-1,0)$. For $1\le q\le L$ on a tail with $L$ resonators, set $b_q:=a_{m-q}$ on the left and $b_q:=a_{m+q}$ on the right. There is a scalar $A$ such that, for every $1\le q\le L$,
\be\lb{tailloc}
b_q=A\big(z_-(\alpha)^q+z_-(\alpha)^{2L+1-q}\big), \quad |b_q|\le|z_-(\alpha)|^{q-1}|b_1|.
\ee
For $1\le q\le L-1$, the adjacent amplitudes satisfy
\be\lb{taildiff}
|b_{q+1}-b_q|\le(1+|z_-(\alpha)|)|z_-(\alpha)|^{q-1}|b_1|.
\ee
\end{prop}

\begin{proof}
Put $z:=z_-(\alpha)$ and $r:=|z|\in(0,1)$, so $z=-r$. For $L=1$, take $A=b_1/[z(1+z)]$; the denominator is nonzero, and there is no adjacent pair within the tail. For $L\ge2$, the recurrence has the general solution $b_q=Az^q+Bz^{-q}$. To impose the free endpoint, introduce $b_{L+1}:=b_L$. Since $b_{L+1}-b_L=0$, we obtain $B=Az^{2L+1}$ after substituting the general solution, and we then obtain the reflected formula in~\eq{tailloc}. In both cases, $b_1=-Ar(1-r^{2L-1})$. If $b_1=0$, the whole tail vanishes. Otherwise,
\[
\frac{b_q}{b_1}=(-r)^{q-1}\frac{1-r^{2L+1-2q}}{1-r^{2L-1}},\quad1\le q\le L.
\]
Since the last fraction lies in $(0,1]$, the bound in~\eq{tailloc} follows. For $1\le q<L$, we have $b_{q+1}-b_q=A(z-1)(z^q-z^{2L-q})$, so
\[
\frac{|b_{q+1}-b_q|}{|b_1|}=(1+r)r^{q-1}\frac{1-r^{2L-2q}}{1-r^{2L-1}}\le(1+r)r^{q-1}.
\]
Thus~\eq{taildiff} holds. If the tail vanishes, both bounds hold directly.
\end{proof}

\begin{proof}[Proof of the localization part of Theorem~\ref{thm:threshold}]
For $0<\tau<\tau_*$, we have $\alpha_N(\tau)>4$ by~\eq{thres}. To normalize its eigenvector at the defect, we first show that $a_m\ne0$. If $a_m=0$, the retained rows of~\eq{perm} would make $\alpha_N(\tau)$ an eigenvalue of $\bL_{N,m}^{\mathrm{re}}(\tau)$, whose spectrum lies below $4$ by Proposition~\ref{prop:band}. We may therefore set $a_m=1$.

The eigenvector decays along each finite tail by Proposition~\ref{prop:tail}. To compute the decay length as $\tau\to\tau_*^-$, put $s:=\alpha_N(\tau)-4\to0^+$. By~\eq{zlocal},
\[
-\log|z_-(4+s)|=2\operatorname{arsinh}(\sqrt{s}/2)=\sqrt{s}\big(1+O(s)\big).
\]
We take the reciprocal and substitute $s=c_*(\tau_*-\tau)+O((\tau_*-\tau)^2)$ from~\eq{crossing}. This proves the stated decay-length asymptotic with relative error $O(\tau_*-\tau)$.
\end{proof}

Because $\ell_*$ describes the homogeneous recurrence, it may exceed a finite tail's length. Formula~\eq{tailloc} still applies, with $N$ fixed. In Section~S3.3 of the supplementary material, we also derive the limiting tail profile at $\tau_*$.

We next bound the amplitudes at the defect neighbors and show that $\alpha_N$ strictly decreases throughout the localized regime. We will use the bound for the gap fields and the strict decrease to recover the radius ratio at a known defect site.

\begin{prop}\lb{prop:monotone}
For each fixed $N$ and defect site $m$, the largest reduced eigenvalue satisfies $\alpha_N'(\tau)<0$ on $(0,\tau_*)$. Its real eigenvector normalized by $a_m=1$ satisfies
\[
-\zeta(\tau)|z_-(\alpha_N(\tau))|\le a_j<0,\quad |j-m|=1,\quad 1\le j\le N.
\]
\end{prop}

\begin{proof}
By Theorem~\ref{thm:threshold}, we can normalize a real eigenvector by $a_m=1$. Put $r:=|z_-(\alpha_N)|\in(0,1)$. On a tail with $L$ resonators, we use~\eq{cont} and Lemma~\ref{lem:tdet} to obtain
\[
\frac{h_{L-1}}{h_L}=-r\frac{1-r^{2L-1}}{1-r^{2L+1}}\in[-r,0),\quad
a_{m\pm1}=\frac{\zeta(\tau)h_{L-1}}{h_L+(\zeta(\tau)-1)h_{L-1}}\in[-\zeta(\tau)r,0).
\]
For the second bound, divide its denominator by $h_L$. Since $\zeta(\tau)\le1$ and $h_{L-1}/h_L<0$, the resulting factor $1+(\zeta(\tau)-1)h_{L-1}/h_L$ is at least one. For each existing neighbor, put $d_j:=1-a_j$. Then $0<d_j<2$, so $d_j^2\le2d_j$, while the defect row gives $\zeta(\tau)\sum_{|j-m|=1}d_j=\alpha_N\tau^3$. We differentiate~\eq{redgeig} along the simple eigenpair, use symmetry to cancel the eigenvector derivatives, and apply these bounds:
\[
\alpha_N'(\tau)\bm a^T\bD_m^{(\tau)}\bm a
=\zeta'(\tau)\sum_{|j-m|=1}d_j^2-3\alpha_N\tau^2
\le\alpha_N\tau^2\Big(\frac{2}{1+\tau}-3\Big)<0.
\]
Here $\zeta'(\tau)/\zeta(\tau)=1/(\tau(1+\tau))$, and every sum includes only existing neighbors. Since $\bm a^T\bD_m^{(\tau)}\bm a>0$, we conclude that $\alpha_N'(\tau)<0$ on $(0,\tau_*)$.
\end{proof}

\subsection{Shrinking-defect asymptotics}
We derived~\eq{redgeig} through~\eq{seq} at fixed $\tau>0$ and now let $\tau\to0$ in this reduced problem, with $N,m$ fixed. For an interior defect, the separation $2\tau R+2\varepsilon$ between $D_{m-1}$ and $D_{m+1}$ stays positive as $\varepsilon\to0$ at fixed $\tau$, but not uniformly as $\tau\to0$. A joint limit would therefore require additional capacitance estimates.

Since $\zeta(\tau)=2\tau(1+O(\tau))$, the defect volume weight $\tau^3$ vanishes faster than its adjacent couplings. We first use this difference in orders to isolate the diverging eigenvalue, then analyze the bounded spectrum.

\begin{lemma}\lb{lem:div}
As $\tau\to0$ in the reduced generalized eigenvalue problem~\eq{redgeig}, exactly one eigenvalue diverges:
\[
\alpha_N=\vartheta_m\frac{\zeta(\tau)}{\tau^3}\big(1+O(\tau^3)\big).
\]
Normalize a corresponding eigenvector $\bm a(\tau)$ by $a_m=1$, and index its components in the original resonator order. Then $\max_{j\ne m}|a_j|=O(\tau^3)$. For each existing neighbor,
\[
a_j=-\frac{\tau^3}{\vartheta_m}+o(\tau^3),\quad 1\le j\le N,\quad |j-m|=1.
\]
\end{lemma}

\begin{proof}
We evaluate the Rayleigh quotient at $\bm e_m\in\bR^N$ and obtain $\alpha_N\ge\vartheta_m\zeta(\tau)/\tau^3\to\infty$. For small $\tau$, $\alpha_N>4$; by Proposition~\ref{prop:band}, it is the only eigenvalue above $4$, and $\lambda_{\max}(\bL_{N,m}^{\mathrm{re}}(\tau))<4$. We can therefore use $g(\alpha_N;\tau)=0$ and bound the resolvent by $(\alpha_N-4)^{-1}$:
\[
0\le\alpha_N-\frac{\vartheta_m\zeta(\tau)}{\tau^3}\le\frac{\vartheta_m\zeta(\tau)^2}{\tau^3(\alpha_N-4)}=O(\zeta(\tau)).
\]
After division by $\vartheta_m\zeta(\tau)/\tau^3$, the relative error is $O(\tau^3)$.

We may set $a_m=1$ by Theorem~\ref{thm:threshold}. Because $\|\bL_{N,m}^{\mathrm{re}}(\tau)\|<4$ and $\alpha_N\to\infty$, we can expand the inverse in the retained rows of~\eq{perm}:
\[
\bm a_m^-=\zeta(\tau)\big(\bL_{N,m}^{\mathrm{re}}(\tau)-\alpha_N\bI_{N-1}\big)^{-1}\bm e_m^{\mathrm{ad}}
=-\frac{\zeta(\tau)}{\alpha_N}\bm e_m^{\mathrm{ad}}+O\Big(\frac{\zeta(\tau)}{\alpha_N^2}\Big)=-\frac{\tau^3}{\vartheta_m}\bm e_m^{\mathrm{ad}}+o(\tau^3).
\]
Since $\bm e_m^{\mathrm{ad}}$ selects exactly the existing neighbors in~\eq{taildata}, this proves both amplitude assertions.
\end{proof}

To compare the bounded eigenvalues with the homogeneous subchain spectra, we now eliminate the defect amplitude from~\eq{perm}.

Set $\bK_m:=\bE_m^{\mathrm{ad}}-\frac1{\vartheta_m}\bm e_m^{\mathrm{ad}}(\bm e_m^{\mathrm{ad}})^T$. For an interior defect, $\bm e_m^{\mathrm{ad}}=\bm e_{m-1}+\bm e_m$, so $\bK_m=\frac12(\bm e_{m-1}-\bm e_m)(\bm e_{m-1}-\bm e_m)^T$. At an endpoint, $\vartheta_m=1$ and $\bE_m^{\mathrm{ad}}=\bm e_m^{\mathrm{ad}}(\bm e_m^{\mathrm{ad}})^T$, so $\bK_m=0$.

Whenever $\vartheta_m\zeta(\tau)-\alpha\tau^3\ne0$, we can solve the defect row for $a_m$. We substitute it into the retained rows and obtain
\be\lb{schurexp}
\begin{aligned}
\bS_m(\alpha,\tau)&:=\bL_{N,m}^{\mathrm{ho}}+\zeta(\tau)\bE_m^{\mathrm{ad}}-\alpha\bI_{N-1} -\frac{\zeta(\tau)^2}{\vartheta_m\zeta(\tau)-\alpha\tau^3}\bm e_m^{\mathrm{ad}}(\bm e_m^{\mathrm{ad}})^T\\
&=\bL_{N,m}^{\mathrm{ho}}+\zeta(\tau)\bK_m-\alpha\bI_{N-1}+\bR_\tau(\alpha),
\end{aligned}
\ee
where
\[
\bR_\tau(\alpha):=-\frac{\alpha\tau^3\zeta(\tau)}{\vartheta_m(\vartheta_m\zeta(\tau)-\alpha\tau^3)}\bm e_m^{\mathrm{ad}}(\bm e_m^{\mathrm{ad}})^T.
\]

\begin{corollary}\lb{cor:slim}
Let $\beta_1\le\cdots\le\beta_{N-1}$ be the eigenvalues of the subchain matrix $\bL_{N,m}^{\mathrm{ho}}$ in~\eq{taildata}, counted with multiplicity. For all sufficiently small $\tau$, the bounded reduced eigenvalues satisfy
\be\lb{brate}
\max_{1\le j\le N-1}|\alpha_j(\tau)-\beta_j|\le C\zeta(\tau),
\ee
where $C$ is independent of $\tau$.
\end{corollary}

\begin{proof}
For small $\tau$ and real $|\alpha|\le5$, $\vartheta_m\zeta(\tau)-\alpha\tau^3\ge\vartheta_m\zeta(\tau)/2$. By~\eq{schurexp} and Weyl's inequality, the ordered eigenvalues $s_j(\alpha,\tau)$ of $\bS_m(\alpha,\tau)$ satisfy $|s_j(\alpha,\tau)-(\beta_j-\alpha)|\le C\zeta(\tau)$. Differentiating~\eq{schurexp} gives
\[
\pa_\alpha\bS_m=-\bI_{N-1}-\frac{\zeta(\tau)^2\tau^3}{(\vartheta_m\zeta(\tau)-\alpha\tau^3)^2}\bm e_m^{\mathrm{ad}}(\bm e_m^{\mathrm{ad}})^T\le-\bI_{N-1}
\]
in the quadratic-form sense. Hence each $s_j$ decreases strictly in $\alpha$ and has opposite signs at $\beta_j\mp M\zeta(\tau)$ for any fixed $M>C$ and sufficiently small $\tau$, since $0\le\beta_j<4$.
The zeros are ordered because the $s_j$ are ordered and strictly decreasing. Schur elimination preserves the kernel dimension, so a shared zero would contradict Proposition~\ref{prop:spec}. Since $\alpha_N$ diverges by Lemma~\ref{lem:div}, the $N-1$ distinct zeros give~\eq{brate}, even when some $\beta_j$ coincide.
\end{proof}

For an interior defect, $\bL_{N,m}^{\mathrm{ho}}$ has a two-dimensional kernel because each subchain supports a constant vector. One combination remains an exact zero mode; the other produces $\alpha_2\to0$. To compute its coefficient and eigenvector, we project the Schur equation onto these two directions.

Put $L:=m-1$ and $K:=N-m$. In $\bR^{N-1}=\bR^L\oplus\bR^K$, set $\cH_0:=\ker\bL_{N,m}^{\mathrm{ho}}$ and $\cH_1:=\cH_0^\perp$, so $\dim\cH_0=2$ and $\dim\cH_1=N-3$. We project onto these spaces with
\[
\bP_0:=\diag(L^{-1}\bm1_L\bm1_L^T,K^{-1}\bm1_K\bm1_K^T),\quad
\bP_1:=\bI_{N-1}-\bP_0.
\]
To represent $\zeta(\tau)\bK_m$ on $\cH_0$, choose the normalized constant vectors
\[
\bm u_L:=L^{-1/2}(\bm1_L,\bm0_K)^T,\quad \bm u_R:=K^{-1/2}(\bm0_L,\bm1_K)^T
\]
as an orthonormal basis of $\cH_0$.

Since $\bK_m=\frac12(\bm e_{m-1}-\bm e_m)(\bm e_{m-1}-\bm e_m)^T$, we represent $\bP_0\bK_m\bP_0$ in the basis $(\bm u_L,\bm u_R)$ by
\be\lb{kcomp}
\frac12\begin{pmatrix}L^{-1}&-(LK)^{-1/2}\\-(LK)^{-1/2}&K^{-1}\end{pmatrix}.
\ee
This $2\times2$ rank-one matrix has eigenvalues $0$ and $(L+K)/(2LK)$.

\begin{lemma}\lb{lem:soft}
Assume the defect is interior. Exactly one positive reduced eigenvalue tends to zero, and as $\tau\to0$ it satisfies
\[
\alpha_2=\frac{L+K}{2LK}\zeta(\tau)\big(1+O(\zeta(\tau))\big).
\]
Put $A_L:=K/\max\{L,K\}$ and $A_R:=-L/\max\{L,K\}$. Normalize the eigenvector by $\max_j|a_j|=1$ and $a_{m-1}>0$. Its components satisfy
\[
a_j=\begin{cases}
A_L+o(1),&1\le j<m,\\
A_R+o(1),&m<j\le N,
\end{cases}\quad a_m=\frac{A_L+A_R}{2}+o(1).
\]
\end{lemma}

\begin{proof}
By Corollary~\ref{cor:slim}, exactly one positive eigenvalue tends to zero, and $0<\alpha_2\le C\zeta(\tau)$. For small $\tau$, we therefore have $2\zeta(\tau)-\alpha_2\tau^3\ge\zeta(\tau)$ and, by~\eq{schurexp}, $\|\bR_\tau(\alpha_2)\|\le C\tau^3$.

To control the component away from $\cH_0$, first normalize $\|\bm a_m^-\|=1$. This is possible because $\bm a_m^-=0$ would force $a_m=0$ in~\eq{perm}. Write $\bm a_{m,i}^-:=\bP_i\bm a_m^-$, $i=0,1$. From the defect row,
\[
a_m=\frac{\zeta(\tau)(\bm e_m^{\mathrm{ad}})^T\bm a_m^-}{2\zeta(\tau)-\alpha_2\tau^3}=O(1).
\]
If $\cH_1\ne\{0\}$, $\bL_{N,m}^{\mathrm{ho}}$ has a fixed positive spectral gap on $\cH_1$. Since $\alpha_2=O(\zeta(\tau))$ and the perturbation in~\eq{schurexp} is $O(\zeta(\tau))$, we can invert the complementary block uniformly for small $\tau$. We solve its equation and use $\|\bm a_m^-\|=1$ to obtain
\[
\|\bm a_{m,1}^-\|=O(\zeta(\tau)),\quad \|\bm a_{m,0}^-\|=1+O(\zeta(\tau)^2).
\]
For $N=3,m=2$, we have $\cH_1=\{0\}$, so these bounds hold with $\bm a_{m,1}^-=0$.

We next isolate the nonzero eigenvalue of~\eq{kcomp}. Since $\tau^3 a_m+\mathbf 1_{N-1}^{T}\boldsymbol a_m^-=0$ by weighted orthogonality and $a_m=O(1)$, $\boldsymbol a_{m,0}^-$ has only an $O(\tau^3)$ component along the global constant direction $(\mathbf 1_L,\mathbf 1_K)$. Its other component lies in the positive eigenspace of~\eq{kcomp}. We project the Schur equation onto $\cH_0$ and take the inner product with $\bm a_{m,0}^-$:
\[
\alpha_2\|\bm a_{m,0}^-\|^2
=\zeta(\tau)\frac{L+K}{2LK}\|\bm a_{m,0}^-\|^2+O(\zeta(\tau)^2+\tau^3).
\]
Since $\|\bm a_{m,0}^-\|^2=1+O(\zeta(\tau)^2)$ and $\zeta(\tau)=2\tau/(1+\tau)$, we can divide by $\|\bm a_{m,0}^-\|^2$ and obtain the claimed relative error $O(\zeta(\tau))$.

Because $\bm a_{m,1}^-\to0$, every subsequential limit of $\bm a_m^-$ is constant on each subchain. If these constants are $c_L,c_R$, the chosen norm, weighted orthogonality, and defect row imply
\[
Lc_L^2+Kc_R^2=1,\quad Lc_L+Kc_R=0,\quad a_m\to\frac{c_L+c_R}{2}.
\]
Thus $|a_{m-1}|\to\sqrt{K/[L(L+K)]}>0$. We choose $a_{m-1}>0$ to fix the limiting sign, then rescale the eigenvector so that its largest absolute component is one. We obtain $A_L,A_R$ and the stated defect amplitude.
\end{proof}

Under the normalization in Lemma~\ref{lem:soft}, the soft mode approaches a constant on each homogeneous subchain, so all differences within either subchain tend to zero. From the left neighbor to the defect, the jump tends to $(A_R-A_L)/2$; from the right neighbor to the defect, it tends to $(A_L-A_R)/2$. Thus both defect-adjacent jumps have nonzero limits although $\zeta(\tau)\to0$.

\begin{proof}[Proof of Theorem~\ref{thm:shrink}]
We obtain~\eq{divbr} from Lemma~\ref{lem:div},~\eq{softbr} from Lemma~\ref{lem:soft}, and the bounded limiting spectrum from Corollary~\ref{cor:slim}. At an endpoint, deleting the defect leaves one homogeneous subchain with a simple zero eigenvalue. Since the full problem already has an exact zero eigenvalue and the bounded eigenvalues converge with multiplicity, no additional positive eigenvalue can tend to zero.
\end{proof}

\section{Resonant fields and narrow-gap gradients}\lb{sec:rec}
We recover the static field at fixed geometry, then use uniform gap estimates to take the narrow-gap limit.

\subsection{Resonant modes for fixed geometry}
At fixed $\varepsilon>0$, choose a simple eigenvalue $\lambda>0$ of~\eq{capgeig} with a real eigenvector $\bm a$ of unit $\bD_m^{(\tau)}$-norm. Let $(\phi_\delta,\psi_\delta)\in X\times X$ be the densities of the corresponding mode in Proposition~\ref{prop:red}. To pass to gap gradients, we need one-sided derivative convergence up to the interfaces, beyond the leading mode approximation in~\cite[Corollary~2.9]{ADH24}.

For every $\varphi\in X$, let $\mathcal{Q}\varphi:=(\int_{\pa D_i}\varphi\,\rd\sigma)_{i=1}^N$ be its complex componentwise charge vector. We use the decomposition in~\cite[Propositions~13--14]{FA22} and define the static kernel coefficients by
\[
\rho_{\bm b}:=\sum_{i=1}^N b_i\psi_i,\quad \mathcal{Q}\rho_{\bm b}=-\bC\bm b,\quad \bm b_\delta:=-\bC^{-1}\mathcal{Q}\phi_\delta.
\]
For every sufficiently small $\delta$, both $\bm b_\delta$ and its weighted inner product with $\bm a$ are nonzero, as we show in the proof of Lemma~\ref{lem:fmode}. We can therefore fix the amplitude and phase by
\be\lb{dnorm}
\bm b_\delta^*\bD_m^{(\tau)}\bm b_\delta=1,\quad \bm b_\delta^*\bD_m^{(\tau)}\bm a>0,
\ee
where the second quantity is required to be real and positive.

By~\eq{statfield}, the static single-layer potential $\cS_D[\rho_{\bm a}]$ equals $\Phi_{\bm a}$ in $\Om$ and the constant $a_i$ throughout each $D_i$.

For each $i$, choose bounded smooth one-sided neighborhoods $U_i^-\subset D_i$ and $U_i^+\subset\Om$ that adjoin $\pa D_i$ and meet no other boundary component.

\begin{lemma}\lb{lem:fmode}
Let $u$ be the corresponding resonant mode normalized by~\eq{dnorm}. For every fixed $0<\sigma<1$ and every pair of compact sets $K_i^\pm\subset\overline{U_i^\pm}$, as $\delta\to0$ we have
\[
\|u-a_i\|_{C^{1,\sigma}(K_i^-)}\to0,\quad\|u-\Phi_{\bm a}\|_{C^{1,\sigma}(K_i^+)}\to0.
\]
The interface norms are one-sided. On compact subsets of $D_i$ and $\Om$, the modes also converge in these norms.
\end{lemma}

\begin{proof}[Proof sketch]
Inverting on the zero-charge complement controls the density correction. The charge equation then gives $\|(\bC-\lambda\bV)\bm b_\delta\|=o_\varepsilon(1)\|\bm b_\delta\|$. Simplicity of $\lambda$ and~\eq{dnorm} then give $\bm b_\delta\to\bm a$. One-sided layer-potential estimates give the stated convergence. See Section~S4.1 of the supplementary material for the details.
\end{proof}

\subsection{Pointwise estimates in narrow regions}
Because the modes converge at fixed geometry as stated above, we can study the narrow-gap limit through the static potentials. We compare each static potential with the linear gap profile, using the graphs $h_j^\pm$, the width $d_{j,j+1}$, and the function $\bar v_j$ from Subsection~\ref{sec:gap}.

\begin{lemma}\lb{lem:nint}
Consider the contact chart for a gap $(D_j,D_{j+1})$, with enlargement $|x'|<2r_0$ separated from the remaining resonators by at least $cR$. Let the adjacent radius ratios satisfy $r_j,r_{j+1}\in[r_-,r_+]$, where $0<r_-\le r_+<\infty$. Let $\bm a=(a_1,\ldots,a_N)^T$ be a complex coefficient vector. In $\Om_{j,j+1}^{r_0}$, the static field~\eq{statfield} satisfies
\[
\Phi_{\bm a}=a_j+(a_{j+1}-a_j)\bar v_j+w_{\bm a,j},\quad
\|\nabla w_{\bm a,j}\|_{L^\infty(\Om_{j,j+1}^{r_0})}\le\frac{C}{R}\|\bm a\|_\infty.
\]
Moreover,
\[
\Big|\sup_{\Om_{j,j+1}^{r_0}}|\nabla\Phi_{\bm a}|-\frac{|a_{j+1}-a_j|}{\varepsilon}\Big|
\le\frac{C}{R}\|\bm a\|_\infty.
\]
The constant is independent of $\varepsilon$ and uniform for radius ratios in $[r_-,r_+]$ under the stated bounds on the normalized charts and separation from the other resonators.
\end{lemma}

\begin{proof}[Proof sketch]
The remainder has zero trace on both faces and satisfies a Poisson equation with source bounded by $C\|\bm a\|_\infty/(R d_{j,j+1})$. Local energy iteration and rescaling by the gap width give the uniform gradient bound. The linear profile has gradient magnitude $\varepsilon^{-1}$ on the contact axis and at most $\varepsilon^{-1}+C/R$ throughout the gap. These two bounds give the supremum estimate. See Section~S4.2 of the supplementary material for the details.
\end{proof}

For $\varepsilon$-dependent amplitudes with $a_{j+1}\ne a_j$, the gradient remainder is of lower order when $(\varepsilon/R)\|\bm a\|_\infty=o(|a_{j+1}-a_j|)$ as $\varepsilon\to0$. Equal neighboring limiting amplitudes give a zero scaled limit in~\eq{moderec}, while an unscaled bound requires control of the jump relative to $\varepsilon$.

\subsection{Sequential field limits and gradient localization}
\begin{proof}[Proof of Theorem~\ref{thm:field}]
Fix $j\in\{2,\ldots,N\}$ and the eigenvector $\bm a^{(j)}$ in the statement. For sufficiently small $\varepsilon$, choose $\bm a_j^\varepsilon$ as in Corollary~\ref{cor:srate}, with unit $\bD_m^{(\tau)}$-norm and positive weighted overlap with $\bm a^{(j)}$. Then $\bm a_j^\varepsilon\to\bm a^{(j)}$ as $\varepsilon\to0$. Normalize $u_j$ by~\eq{dnorm}, using $\bm a_j^\varepsilon$ as the reference eigenvector.

\textbf{Step 1: Sequential gradient limit.} For each fixed $\varepsilon>0$, we set $\Phi:=\Phi_{\bm a_j^\varepsilon}$ and first compare $u_j$ with this static field. Since the closed gap chart has a finite cover by one-sided boundary neighborhoods and compact exterior sets inside the larger chart, Lemma~\ref{lem:fmode} applies on this cover. Hence
\[
\lim_{\delta\to0}\varepsilon\|\nabla(u_j-\Phi)\|_{L^\infty(\Om_{\ell,\ell+1}^{r_0})}=0.
\]
We then apply Lemma~\ref{lem:nint} to obtain
\[
\Big|\varepsilon\sup_{\Om_{\ell,\ell+1}^{r_0}}|\nabla\Phi|
-|(\bm a_j^\varepsilon)_{\ell+1}-(\bm a_j^\varepsilon)_\ell|\Big|
\le C\frac{\varepsilon}{R}\|\bm a_j^\varepsilon\|_\infty.
\]
We next let $\varepsilon\to0$ with $\tau$ fixed. Since $\bm a_j^\varepsilon$ has unit $\bD_m^{(\tau)}$-norm, $\|\bm a_j^\varepsilon\|_\infty\le\tau^{-3/2}$ independently of $\varepsilon$, and the right-hand side tends to zero. We therefore obtain~\eq{moderec} by taking the limits in the order~\eq{seq} and using $\bm a_j^\varepsilon\to\bm a^{(j)}$.

\textbf{Step 2: Gradient localization of the largest-eigenvalue mode.} Fix $0<\tau<\tau_*$. On a tail of length $L$, let $b_q$ denote the reduced amplitude at the $q$th resonator counted outward from the defect, divided by $a_m^{(N)}$, for $1\le q\le L$. Put $r:=|z_-(\alpha_N(\tau))|\in(0,1)$. On each tail, Proposition~\ref{prop:monotone} gives $-\zeta(\tau)r\le b_1<0$. Thus $|1-b_1|\le1+\zeta(\tau)r\le2$ across the gap adjacent to the defect. Between successive resonators within a tail ($2\le q\le L$),~\eq{taildiff} gives
\[
|b_q-b_{q-1}|\le(1+r)r^{q-2}|b_1|\le(1+r)\zeta(\tau)r^{q-1}\le2r^{q-1}.
\]
A tail with one resonator has no internal gap, so only the first bound is needed. Applying~\eq{moderec} to $\widetilde u_N$ proves~\eq{phystail}.
\end{proof}

We finally let the defect shrink. Since Lemma~\ref{lem:div} controls every amplitude away from the defect, we can evaluate each gap-gradient limit as follows.

\begin{corollary}\lb{cor:grad}
For the rescaled largest-eigenvalue mode in Theorem~\ref{thm:field}, let $q$ denote the outward gap index defined there. Choose $r_0=r_0(\tau)$ as in Subsection~\ref{sec:setup}, independently of $\varepsilon$ at each fixed $\tau>0$. For every $1\le\ell<N$,
\[
\lim_{\tau\to0}\lim_{\varepsilon\to0}\varepsilon\lim_{\delta\to0}
\sup_{\Om_{\ell,\ell+1}^{r_0}}|\nabla\widetilde u_N|
=\begin{cases}
1, & q=1,\\
0, & q\ge2.
\end{cases}
\]
\end{corollary}

\begin{proof}
Lemma~\ref{lem:div} gives $a_j^{(N)}/a_m^{(N)}=O(\tau^3)$ for every $j\ne m$. The normalized amplitude difference $|a_{\ell+1}^{(N)}-a_\ell^{(N)}|/|a_m^{(N)}|$ therefore tends to one in a defect-adjacent gap and to zero in every other gap. Applying~\eq{moderec} to $\widetilde u_N$ before letting $\tau\to0$ proves both cases.
\end{proof}

\section{Conclusions}
The exact characteristic determinant determines the finite-chain localization threshold and underlies the shrinking-defect asymptotics. The absolute jumps in reduced amplitudes determine the sequential scaled gap gradients and their concentration next to the shrinking defect for the largest-eigenvalue mode. With several defects, localized modes may interact and split in frequency depending on defect separation and position. Joint limits in $\tau$, $\varepsilon$, and $\delta$ require uniform capacitance and gap-field estimates and may yield different coupling regimes. For nonspherical resonators, the couplings must be derived from the local contact geometry. Relating resonance frequencies and gap fields to driven scattering measurements could help recover the defect size and position.

\bigskip\noindent{\bf Acknowledgments.}
H. Li was partially supported by the National Natural Science Foundation of China (Grant No.~12471191) and the Beijing Natural Science Foundation (Grant No.~1242006).

\end{document}


\title[Supplementary material]
{Supplementary material for\\
Subwavelength resonances and localization in finite resonator chains with a size defect}
\author[Y.~Gong, H.~Li, and J.~Zhang]{Yihao Gong, Haigang Li, and Junhua Zhang}
\maketitle

We use the notation and assumptions of the main paper. Sections~S1 and~S2 give details of the capacitance reduction, Section~S3 gives additional spectral results and comparisons, and Section~S4 supplies the field estimates. Numbers prefixed by S refer to this supplement.

\section{Layer-potential formulation and fixed-geometry reduction}\lb{si:layer}

Fix $\varepsilon>0$ and $\tau\in(0,1]$ as in Subsection~\ref{M-sec:setup}, with the layer operators, spaces $X,Y$, and trace conventions of Section~2.1 of the main paper. The kernel satisfies $(\Delta_x+k^2)G^k(x,y)=\delta_y$, the Dirac mass at $y$.

Using the jump relations~\cite[Chapters~1--2]{AKL09}, the transmission conditions give~\cite[Section~2A]{ADH24}
\be\lb{si:layersys}
\begin{pmatrix}
\cS_D^{k_b}&-\cS_D^k\\
-\frac12\cI+\cK_D^{k_b,*}&-\delta(\frac12\cI+\cK_D^{k,*})
\end{pmatrix}\binom{\phi}{\psi}=0
\ee
in $Y\times X$, with interior and exterior densities $\phi,\psi\in X$, respectively. At $k=0$ the trace $\cS_D:X\to Y$ is invertible~\cite[Theorems~1--2 and Section~5]{Co88}. Ellipticity lifts its inverse to $H^{s+1}(\pa D)\to H^s(\pa D)$ for $s\ge-1/2$, since the diagonal blocks of $\cS_D$ act on spherical harmonics of degree $\ell\ge0$ with multipliers $-R_i/(2\ell+1)$ and the off-diagonal kernels are smooth.

The charge representation of $\bC$ and the static field representation are given in~\eq{M-statfield}.

For each fixed $s\ge-1/2$ and small complex $k$, the kernel expansion and smooth-boundary mapping estimates of~\cite[Proposition~12]{FA22} and~\cite[Example~2.4 and Theorem~2.6]{CL12} give
\[
\|\cS_D^k-\cS_D\|_{H^s\to H^{s+1}}=O_{\varepsilon,s}(|k|),\quad
\|\cK_D^{k,*}-\cK_D^*\|_{H^s\to H^s}=O_{\varepsilon,s}(|k|^2).
\]
Here the Sobolev spaces are over $\pa D$; the linear term of $G^k$ has zero normal derivative. Perturbation preserves the invertibility of $\cS_D^k:H^s(\pa D)\to H^{s+1}(\pa D)$ for small $k$.

\subsection{Proof of Proposition~\ref{M-prop:red}}\lb{si:fixed-reduction}

\begin{proof}
The generalized capacitance matrix is $\delta v_b^2\bV^{-1}\bC$~\cite[Definition~2.6]{ADH24}, so Theorem~2.7 there gives the stated expansion with the positive square root. For uniqueness, the rescaled matrix~\cite[(154)]{FA22} equals $\lambda\bV-\bC$ at $\om=0$; the proof of Proposition~18 there gives an analytic branch $\Lambda(\om)$ with $\Lambda(0)=\lambda>0$.

Writing $\om=\sqrt\delta\,\xi$, we solve $\xi^2-v_b^2\Lambda(\sqrt\delta\,\xi)=0$. At $(\delta,\xi)=(0,v_b\sqrt\lambda)$, the left-hand side has $\xi$-derivative $2v_b\sqrt\lambda\ne0$. The implicit function argument~\cite[proof of Proposition~20]{FA22}, with parameter $\sqrt\delta$, gives the unique branch $\xi=v_b\sqrt\lambda+O(\sqrt\delta)$ at fixed geometry, consistently with~\eq{M-seq}.
\end{proof}

\section{Capacitance and spectral approximation estimates}\lb{si:cap}

\subsection{Proof of Lemma~\ref{M-lem:neck}}\lb{si:neck-proof}
\begin{proof}
Put $t:=|x'|$ and $a_{rs}:=(r+s)/(2Rrs)$. Expanding the spherical graphs gives $d_{rs}(x')=\varepsilon+a_{rs}t^2+E_{rs}(t)$ with
\be\lb{si:neckwidth}
0\le E_{rs}(t)\le Ct^4/R^3,\quad\text{and}\quad
C^{-1}(\varepsilon+|x'|^2/R)\le d_{rs}(x')\le C(\varepsilon+|x'|^2/R),
\ee
uniformly for $r,s\in[r_-,r_+]$ and $t<2r_0$. The quadratic reciprocal integrates to $\pi a_{rs}^{-1}\log(1+a_{rs}r_0^2/\varepsilon)$ on $|x'|<r_0$. By~\eq{si:neckwidth} the reciprocal error is bounded by $Ct^4/[R^3(\varepsilon+t^2/R)^2]\le C/R$, giving an $O(R)$ integral error. Since $a_{rs}r_0^2/R$ is bounded above and away from zero, we obtain $J_{rs}(\varepsilon)=\pi a_{rs}^{-1}|\log(\varepsilon/R)|+O(R)$ with $\pi a_{rs}^{-1}=2\pi Rrs/(r+s)$.

Integrating $\pa_{x^1}u$ between the traces and using Cauchy--Schwarz gives~\eq{M-fiblo}. The linear profile has longitudinal energy $|A-B|^2J_{rs}$, and differentiation as in~\cite[(2.3)--(2.5)]{LZ23} gives $|\nabla_{x'}\bar v_{A,B}|\le C|A-B|t/(Rd_{rs})$. Integrating across the gap bounds its transverse energy by $C|A-B|^2R^{-2}\int_{|x'|<r_0}t^2/d_{rs}\,\rd x'\le CR|A-B|^2$.
\end{proof}

\subsection{Proof of Proposition~\ref{M-prop:cap}}\lb{si:cap-proof}

\begin{proof}
Fix $N,R$ and let $\tau$ range over the compact set of Proposition~\ref{M-prop:cap}, with its uniform chart and separation bounds. Restricting the energy of $\Phi_{\bm a}$ to the disjoint gaps gives the lower bound~\eq{M-fiblo}.

For the upper bound, use $a_j+(a_{j+1}-a_j)\bar v_j$ in the enlarged gap chart. In $r_0\le |x'|<2r_0$ the width is at least $cR$, so the profile and its gradient are bounded by $C\|\bm a\|$ and $C\|\bm a\|/R$, respectively. Extend the constant traces through the separated collars with the same bounds and join them to the profiles using cutoffs of gradient $O(R^{-1})$. This gives a compactly supported $W_{\bm a}\in H^1(\Om)$, linear in $\bm a$, with the prescribed traces, the same gap profiles, and support of volume $O(R^3)$. Hence
\[
\int_{\Om\setminus\bigcup_{j=1}^{N-1}\Om_{j,j+1}^{r_0}}|\nabla W_{\bm a}|^2\,\rd x\le CR\|\bm a\|^2.
\]

Since the traces agree and the far-field term vanishes by decay, the exterior Dirichlet principle~\cite[(1.20)--(1.21) and Lemma~2.2]{LZ23} gives $\bm a^T\bC\bm a\le\int_\Om|\nabla W_{\bm a}|^2\,\rd x$. Lemma~\ref{M-lem:neck} then gives
\[
\Big|\bm a^T\bC\bm a-\sum_{j=1}^{N-1}|a_{j+1}-a_j|^2J_{r_j,r_{j+1}}\Big|\le CR\|\bm a\|^2.
\]
Substituting the gap-energy expansion from Section~2.2 of the main paper and taking the supremum over real unit vectors gives~\eq{M-capexp}.
\end{proof}

\subsection{Proof of Corollary~\ref{M-cor:srate}}\lb{si:spectral-proof}

\begin{proof}
Multiply the matrix difference in~\eq{M-capexp} on both sides by $\bV^{-1/2}$ and use~\eq{M-volmat} and the normalization of $\widehat{\bH}_\varepsilon$ to obtain~\eq{M-Hrate}. Let $\eta_\varepsilon:=\|\widehat{\bH}_\varepsilon-\widehat{\bH}(\tau)\|$ and let $\gamma:=\min_{1\le j<N}(\alpha_{j+1}-\alpha_j)>0$ be the smallest reduced spectral gap. Weyl's inequality gives~\eq{M-srate} and simplicity when $\eta_\varepsilon<\gamma/3$.

For the Euclidean unit vectors $\bm c_j^\varepsilon:=(\bD_m^{(\tau)})^{1/2}\bm a_j^\varepsilon$ and $\bm c^{(j)}:=(\bD_m^{(\tau)})^{1/2}\bm a^{(j)}$, the perturbed eigenvalue equation gives
\[
\|(\widehat{\bH}(\tau)-\alpha_j\bI)\bm c_j^\varepsilon\|
\le|\mu_{j,\varepsilon}-\alpha_j|+\eta_\varepsilon\le 2\eta_\varepsilon.
\]
The reduced spectral gap bounds the component orthogonal to $\bm c^{(j)}$ by $2\eta_\varepsilon/\gamma$. With positive overlap this bounds $\|\bm c_j^\varepsilon-\bm c^{(j)}\|$ by $C\eta_\varepsilon/\gamma$. Multiplying by $(\bD_m^{(\tau)})^{-1/2}$ gives the asserted estimate. This inverse transform is uniformly bounded and $\gamma$ stays bounded away from zero on the stated compact intervals.
\end{proof}

\section{Additional spectral consequences and consistency checks}\lb{si:additional}

\subsection{Stability of the band crossing}\lb{si:band-stability}

\begin{corollary}\lb{si:cor-band-stability}
Fix $N,m,R$ and a compact interval $I\Subset(0,1)$ containing $\tau_*$ in its interior, and assume uniform geometric bounds. Uniformly for $\tau\in I$, the largest normalized capacitance eigenvalue $\mu_{N,\varepsilon}$ from Corollary~\ref{M-cor:srate} satisfies
\be\lb{si:bandstab}
\mu_{N,\varepsilon}(\tau)-4=c_*(\tau_*-\tau)+O\big((\tau-\tau_*)^2+|\log(\varepsilon/R)|^{-1}\big).
\ee
There exist constants $K_0,\eta,\varepsilon_0>0$ such that $[\tau_*-\eta,\tau_*+\eta]\subset I$ and, whenever
\[
0<\varepsilon<\varepsilon_0,\quad K_0|\log(\varepsilon/R)|^{-1}\le|\tau-\tau_*|\le\eta,
\]
we have $\mu_{N,\varepsilon}(\tau)>4$ for $\tau<\tau_*$ and $\mu_{N,\varepsilon}(\tau)<4$ for $\tau>\tau_*$.
\end{corollary}

\begin{proof}
Combine the Taylor expansion~\eq{M-crossing}, whose quadratic remainder is uniform by analyticity of $\alpha_N$ on $I$, with Corollary~\ref{M-cor:srate}. In the stated range, small $\eta$ and large $K_0$ make each error at most $c_*|\tau-\tau_*|/4$, so the linear term determines the sign. Choose $\varepsilon_0$ so that $K_0|\log(\varepsilon/R)|^{-1}\le\eta$ for $0<\varepsilon<\varepsilon_0$.
\end{proof}

\subsection{Algebraic extension to enlarged defects}\lb{si:large-defects}

\begin{remark}\lb{si:remark-large-defects}
For the algebraic extension $\tau\ge1$, put $w:=\zeta(\tau)\in[1,2)$. In the order from the defect outward, a retained block is $\bL_{L,0}+w\bm e_1\bm e_1^T$. For nonzero $\bm b\in\bR^L$,
\[
\bm b^T(4\bI_L-\bL_{L,0}-w\bm e_1\bm e_1^T)\bm b
=\sum_{q=1}^{L-1}(b_q+b_{q+1})^2+2b_L^2+(2-w)b_1^2>0.
\]
Indeed, vanishing of the sum and the $b_L$ term forces $\bm b=0$. The Schur argument of Proposition~\ref{M-prop:band} applies on $[4,\infty)$. The evaluation~\eq{M-phisign} remains valid since $4L-(2L-1)w>0$ for each tail. Since $H$ decreases, $g(4;\tau)=2\tau H(\tau)\le2\tau H(1)<0$. Hence the extended reduced problem has no eigenvalue at or above $4$.
\end{remark}

\subsection{The critical finite-tail profile}\lb{si:critical-tail}

\begin{remark}\lb{si:remark-critical-tail}
For $\tau_*<\tau\le1$, the unit-modulus conjugate roots give oscillatory tails for $\alpha_N(\tau)$. At $\tau=\tau_*$ the repeated root $-1$ gives $b_q=(-1)^q(A+Bq)$ subject to the free-endpoint condition.

The shifted retained block at $\alpha_N=4$ is invertible by Proposition~\ref{M-prop:band} and remains so nearby. Hence $a_m\ne0$ and we set $a_m=1$. On a tail of length $L$, Lemma~\ref{M-lem:tdet} then gives $b_1=\zeta(\tau)h_{L-1}(\alpha_N)/q_L(\alpha_N;\tau)$, nonzero near the crossing since $h_{L-1}(4)\ne0$. Backward recursion from the free endpoint and $h_s(4)=(-1)^s(2s+1)$ give, as $\tau\to\tau_*$ from either side,
\[
\frac{b_q}{b_1}=\frac{h_{L-q}(\alpha_N)}{h_{L-1}(\alpha_N)}
\longrightarrow(-1)^{q-1}\frac{2L+1-2q}{2L-1},\quad 1\le q\le L.
\]
\end{remark}

\subsection{Homogeneous-chain and two-sphere checks}\lb{si:checks}

At $\tau=1$ and $N>2$, the homogeneous spectrum and the leading capacity term in~\eq{M-cmean} agree with~\cite[Theorem~1.1 and Lemma~3.5]{LZ25}. For $N=2$, the reduced determinant gives $\alpha_2(\tau)=\zeta(\tau)(1+\tau^{-3})$, which yields the logarithmic frequency coefficient of~\cite[Theorem~4.4]{ADY20}. At equal radii, $\operatorname{Cap}(D)\to8\pi R\log2$ by~\cite[(3.1)--(3.4)]{ADY20}. Together with~\eq{M-cmean}, this gives
\[
\alpha_2(1)=2,\quad
\lim_{\varepsilon\to0}\lim_{\delta\to0}\frac{\om_1^2}{\delta v_b^2}=\frac{3\log2}{R^2}.
\]

With the defect first and $a_1=1$, the positive-eigenvalue vector is $(1,-\tau^3)^T$. For this limiting amplitude,~\eq{M-moderec} gives scaled gap gradient $1+\tau^3$, with the $\varepsilon^{-1}$ rate of~\cite[Theorem~5.2]{ADY20}. These comparisons use the sequential order~\eq{M-seq}; the cited PDE results use coupled limits.

\section{Technical proofs for resonant fields}\lb{si:fields}

\subsection{Convergence of resonant modes at fixed geometry}\lb{si:mode-proof}

\begin{proof}[Proof of Lemma~\ref{M-lem:fmode}]
Fix the geometry and eigenpair $(\lambda,\bm a)$ of Lemma~\ref{M-lem:fmode}, with the charges $\mathcal Qf$ and densities $\rho_{\bm b}$ of Section~5.1 of the main paper. The projection onto $\operatorname{span}\{\psi_i\}$ along $\ker\mathcal Q$ is~\cite[Propositions~13--14 and (125)--(127)]{FA22}
\[
\mathcal H f:=\rho_{-\bC^{-1}\mathcal Q f}.
\]
Indeed, $\mathcal Q\rho_{\bm b}=-\bC\bm b$, so $\mathcal H\phi_\delta=\rho_{\bm b_\delta}$ with $\bm b_\delta=-\bC^{-1}\mathcal Q\phi_\delta$. Since $\ker(-\frac12\cI+\cK_D^*)=\operatorname{span}\{\psi_1,\ldots,\psi_N\}$, the operator $\mathcal L_0:=-\cI/2+\cK_D^*+\mathcal H$ acts as the identity on this span and agrees with $-\cI/2+\cK_D^*$ on $\ker\mathcal Q$. Its $L^2$ invertibility extends to $H^s(\pa D)$, $s\ge-1/2$: the one-order gain of $\cK_D^*$~\cite[Example~2.4 and Theorem~2.6]{CL12} and the smooth range of $\mathcal H$ give Fredholm index zero and a smooth, hence trivial, kernel.

Eliminating $\psi_\delta=(\cS_D^k)^{-1}\cS_D^{k_b}\phi_\delta$ from~\eq{si:layersys} and adding $\mathcal H$ gives the operator of~\cite[(126)--(130)]{FA22}:
\[
\mathcal G(\om,\delta):=-\frac12\cI+\cK_D^{k_b,*}-\delta\Big(\frac12\cI+\cK_D^{k,*}\Big)(\cS_D^k)^{-1}\cS_D^{k_b}+\mathcal H.
\]
Section~S1 gives $\|\mathcal G-\mathcal L_0\|_{H^s\to H^s}\le C_{\varepsilon,s}(|\om|^2+\delta)$. Thus $\mathcal G$ is invertible for small $\delta$ in each fixed $H^s$, and the inverses agree on smooth data by uniqueness in $H^{-1/2}$. The homogeneous equation~\cite[(136)--(137)]{FA22} gives
\[
\mathcal G\phi_\delta=\mathcal H\phi_\delta=\rho_{\bm b_\delta},
\quad \phi_\delta=\mathcal G^{-1}\rho_{\bm b_\delta}.
\]
If $\bm b_\delta=0$, both densities vanish by this identity and the first row of~\eq{si:layersys}. Since $\mathcal L_0^{-1}\rho_{\bm b_\delta}=\rho_{\bm b_\delta}$ and $|\om|=O_\varepsilon(\sqrt\delta)$, the inverse expansion~\cite[(130)]{FA22} gives
\[
\|\phi_\delta-\rho_{\bm b_\delta}\|_{H^s}\le C_{\varepsilon,s}\delta\|\bm b_\delta\|.
\]
The first row of~\eq{si:layersys} and Section~S1 also give $\|\psi_\delta-\phi_\delta\|_{H^s}\le C_{\varepsilon,s}|\om|\|\phi_\delta\|_{H^s}$.

The charge vector in~\cite[(135)--(138)]{FA22} is $-\bC\bm b_\delta$. Substituting it into the homogeneous system with the expansion in~\cite[Proposition~17, (141)]{FA22} yields
\[
\Big\|\Big(\frac{\om^2}{v_b^2}\bV-\delta\bC\Big)\bm b_\delta\Big\|
\le C_\varepsilon\Big(|\om|^3+\delta|\om|+(|\om|^2+\delta)^2\Big)\|\bm b_\delta\|.
\]
Dividing by $\delta$ and using Proposition~\ref{M-prop:red} gives $\|(\bC-\lambda\bV)\bm b_\delta\|=o_\varepsilon(1)\|\bm b_\delta\|$.

Normalize $\bm b_\delta$ in the $\bD_m^{(\tau)}$-norm, proportional to the $\bV$-norm by~\eq{M-volmat}. The fixed spectral gap on the $\bV$-orthogonal complement of $\bm a$ gives $|\bm b_\delta^*\bD_m^{(\tau)}\bm a|\to1$. The phase in~\eq{M-dnorm} yields $\bm b_\delta\to\bm a$ and both densities converge to $\rho_{\bm a}$ in each fixed $H^s$.

On $U=U_i^\pm$, the kernel expansion and~\cite[Theorem~2.7]{CL12} give
\[
\|\cS_D^k-\cS_D\|_{H^s(\pa D)\to H^{s+3/2}(U)}\le C_{\varepsilon,s,U}|k|.
\]
The potentials thus converge in $H^{s+3/2}(U_i^\pm)$ to $\cS_D[\rho_{\bm a}]$, equal to $\Phi_{\bm a}$ outside and $a_i$ inside. Taking $s>1+\sigma$, Sobolev embedding gives one-sided $C^{1,\sigma}$ convergence; away from the interfaces, the kernels are smooth.
\end{proof}

\subsection{Uniform estimates in a narrow gap}\lb{si:gap-proof}

\begin{proof}[Proof of Lemma~\ref{M-lem:nint}]
Put $A:=\|\bm a\|_\infty$ and write $w:=w_{\bm a,j}=\Phi_{\bm a}-a_j-(a_{j+1}-a_j)\bar v_j$ in $\Om_{j,j+1}^{2r_0}$. The maximum principle for $|\Phi_{\bm a}|^2$ and $0\le\bar v_j\le1$ give $|w|\le2A$. The remainder vanishes on both faces and satisfies $\Delta w=-(a_{j+1}-a_j)\Delta\bar v_j$. Since the slopes of the graphs are $O(|x'|/R)$ and their second derivatives are $O(R^{-1})$, differentiating the profile gives $|\Delta\bar v_j|\le C/(Rd_{j,j+1})$ by~\eq{si:neckwidth}.

Scale by $R$ and write $W(y):=w(Ry)$, $d(y'):=d_{j,j+1}(Ry')/R$ and $\widehat r_0:=r_0/R$. For the scaled graphs $H^\pm(y'):=h_j^\pm(Ry')/R$ and a local radius $t>0$, set
\[
\mathcal U_t(y_0'):=\{y:|y'-y_0'|<t,\ H^-(y')<y^1<H^+(y')\}.
\]
Then $\Delta_yW=F$ with $|F|\le CA/d(y')$, and~\eq{si:neckwidth} gives
\[
C^{-1}(\varepsilon/R+|y'|^2)\le d(y')\le C(\varepsilon/R+|y'|^2).
\]
Testing with $\eta^2\overline W$ and taking real parts gives $\int_{\mathcal U_{3\widehat r_0/2}(0)}|\nabla_yW|^2\,\rd y\le CA^2$. Here $\eta$ is a tangential cutoff supported in $|y'|<2\widehat r_0$ and equal to one on $|y'|\le3\widehat r_0/2$. Integration across the gap gives $\|F\|_{L^1(\mathcal U_{2\widehat r_0}(0))}\le CA$.

For $|y_0'|<\widehat r_0$ and small $d_0:=d(y_0')$, the width satisfies $C^{-1}d_0\le d(y')\le Cd_0$ on $|y'-y_0'|\le c_0\sqrt{d_0}$ with fixed small $c_0$. Set $E(t):=\int_{\mathcal U_t(y_0')}|\nabla_yW|^2\,\rd y$. The zero traces give the Poincar\'e factor $Cd_0^2$, and $\|F\|_{L^2(\mathcal U_s(y_0'))}^2\le CA^2s^2/d_0$. The cutoff estimate in~\cite[proof of Proposition~1.7, (2.7)--(2.10)]{Li20} therefore gives
\[
E(t)\le C\frac{d_0^2}{(s-t)^2}E(s)+CA^2s^2d_0,\quad 0<t<s\le c_0\sqrt{d_0}.
\]
Taking steps $t_i:=2d_0+iMd_0$ with fixed large $M$ yields
\[
E(t_i)\le\theta E(t_{i+1})+CA^2(i+1)^2d_0^3,\quad t_{i+1}\le c_0\sqrt{d_0},\quad 0<\theta<1.
\]
The cylinders stay inside the initial energy region. With at least $c/\sqrt{d_0}$ steps, the cited iteration gives $E(2d_0)\le CA^2d_0^3$.

Rescale by $d_0$ to obtain uniformly separated $C^{1,1}$ faces. Since $d_0\|F\|_{L^\infty(\mathcal U_{2d_0}(y_0'))}\le CA$, the estimate in~\cite[proof of Proposition~1.7, (2.11)]{Li20}, applied to the real and imaginary parts, gives
\[
\|\nabla_yW\|_{L^\infty(\mathcal U_{d_0}(y_0'))}\le C\big(d_0^{-3/2}E(2d_0)^{1/2}+A\big)\le CA.
\]
Where $d_0$ is bounded below, fixed-scale estimates in the enlarged chart give the same bound. Scaling back gives $\|\nabla w\|_{L^\infty(\Om_{j,j+1}^{r_0})}\le CA/R$, independently of $\varepsilon$ under the stated radius, chart, and separation bounds.

The spherical graphs give $|\nabla_{x'}h_j^-+\bar v_j\nabla_{x'}d_{j,j+1}|\le C|x'|/R$. Thus
\[
|\nabla\bar v_j|=\frac{\sqrt{1+|\nabla_{x'}h_j^-+\bar v_j\nabla_{x'}d_{j,j+1}|^2}}{d_{j,j+1}}
\le\varepsilon^{-1}+C/R
\]
by the lower bound $d_{j,j+1}\ge\varepsilon+c|x'|^2/R$. On $x'=0$, the magnitude is $\varepsilon^{-1}$. The triangle inequality in the gap and its reverse on the axis give the supremum estimate, since $|a_{j+1}-a_j|\le2A$.
\end{proof}